\documentclass[10pt]{article}
\usepackage[square,authoryear]{natbib}
\usepackage{marsden_article}
\usepackage{epstopdf}
\usepackage{setspace}
\usepackage{bbm}

\usepackage{amscd}

\usepackage{caption2}

\DeclareGraphicsExtensions{.eps,.ps,.pdf}

\begin{document}

 \newtheorem{thm}{Theorem}[section]
 \newtheorem{cor}[thm]{Corollary}
 \newtheorem{lem}[thm]{Lemma}
 \newtheorem{prop}[thm]{Proposition}
 \newtheorem{defn}[thm]{Definition}
 \newtheorem{rem}[thm]{Remark}
 \numberwithin{equation}{section}


\title{\Large{\textbf{A Geometric Theory of Growth Mechanics}}
\footnote{Dedicated to the memory of Professor James K. Knowles (1931-2009).}
\footnote{To appear in the \emph{Journal of Nonlinear Science}.}}

\author{ Arash Yavari\thanks{School of Civil and Environmental Engineering,
  Georgia Institute of Technology, Atlanta, GA 30332. E-mail: arash.yavari@ce.gatech.edu.}}

\maketitle

\begin{abstract}
In this paper we formulate a geometric theory of the mechanics of
growing solids. Bulk growth is modeled by a material manifold with
an evolving metric. Time dependence of metric represents the evolution of the stress-free (natural) configuration of the body in response to changes
in mass density and ``shape". We show that time dependency of
material metric will affect the energy balance and the entropy
production inequality; both the energy balance and the entropy
production inequality have to be modified. We then obtain the
governing equations covariantly by postulating invariance of
energy balance under time-dependent spatial diffeomorphisms. We use the principle of maximum entropy production in deriving an evolution equation for the material metric. In
the case of isotropic growth, we find those growth distributions
that do not result in residual stresses. We then look at Lagrangian field theory of growing elastic solids. We will use the Lagrange-d'Alembert's principle with Rayleigh's dissipation functions to derive all the governing equations. We make an explicit
connection between our geometric theory and the conventional
multiplicative decomposition of deformation gradient
$\mathbf{F}=\mathbf{F}_e\mathbf{F}_g$ into growth and elastic
parts. We linearize the nonlinear theory and derive a linearized
theory of growth mechanics. Finally, we obtain the stress-free
growth distributions in the linearized theory.
\end{abstract}


\tableofcontents


\section{Introduction}

In classical continuum mechanics, one usually models
mass-conserving bodies. The traditional framework of continuum
mechanics is suitable for many practical applications. However, in
some natural phenomena mass is added or lost in a deformation
process. This is particularly important in biological systems
where growth and remodeling are closely linked to mechanical
loads. In the case of soft tissues, elastic deformations are large
and theory of linear elasticity is not adequate. This has been the
main motivation for the extensive studies of large deformations in
biomechanics in recent years (see
\citet{CowinHegedus1976,Skalak1982,Fung1983,Naumov1994,Hoger1997,Humphrey2003,KlarbringOlssonStalhand2007}
and references therein).

It has been realized for a long time that mechanical forces
directly affect growth and remodeling in biological systems
\citep{Hsu1968}. A continuum theory capable of modeling biological
tissues must take into account changes of mass and the coupling
between growth/remodeling and mechanical stresses. In continuum
mechanics, one starts by postulating that a body is made of a
large number of ``material points" that can be treated as
mathematical points. A material point consists of ``enough" number
of particles (atoms, molecules, cells, etc.) such that it can
represent the mechanical properties of the body, e.g. density.
Material points are then identified with their positions in the
so-called reference configuration. This is called the material
manifold. It is always assumed that the body is macroscopically
stress-free in the material manifold. The material manifold is not
necessarily Euclidean and even not Riemannian; in general,
material manifold is a Riemann-Cartan manifold in the case of
solid bodies with distributed dislocations, for example. It is
relevant to mention that in most of the existing formulations of
finite-strain plasticity, instead of working with a Riemann-Cartan
material manifold, one assumes a multiplicative decomposition of
the deformation gradient into elastic and plastic parts, i.e.
\citep{Eckart1948,Kroner1960,LeeLiu1967,Lee1969}
\begin{equation}
    \mathbf{F}=T\varphi=\mathbf{F}_e\mathbf{F}_p,
\end{equation}
where $\varphi$ is the deformation mapping. This means that
locally the material deforms plastically followed by elastic
deformations to ensure compatibility. In other words, one assumes
that both the material and the ambient space manifolds are
Euclidean and one locally decomposes the total deformation mapping
into incompatible elastic and plastic parts. Motivated by
plasticity, in the case of growth several researchers
\citep{KondaurovNikitin1987,TakamizawaMatsuda1990,Takamizawa1991,Rodriguez1994}
postulated a similar decomposition of $\mathbf{F}$ into elastic
and growth parts, i.e.
\begin{equation}
    \mathbf{F}=\mathbf{F}_e\mathbf{F}_g.
\end{equation}
This has been the fundamental kinematical idea of all the existing
models of growth mechanics to this date (see
\citet{BenAmarGoriely2005,Hoger1997,LubardaHoger2002} and
references therein).

Recently, \cite{OzakinYavari2009} introduced a geometric theory of
thermoelasticity in which thermal strains are buried in a
temperature-dependent Riemannian material manifold. In that theory
a change of temperature leads to a rescaling of the material
metric with a clear physical meaning. In this paper we introduce a
geometric theory for growing bodies using similar ideas. We should
mention that the analogy between growth and thermal distortions
was first realized by Skalak and his coworkers \citep{Skalak1996}.

There are two possibilities in a growth process: (i) the number of
material points is fixed, and (ii) material points are removed or
are added. Note that in a continuum model material points are
assumed to contain several (or a large number of) particles,
cells, etc. Erosion or accretion of cells corresponds to changes
in volume (and the corresponding mass) and shape of the material
body. In our continuum model, similar to many of the earlier
models on bulk growth, we assume that the number of material
points is fixed. This means that we work with a fixed set
$\mathcal{B}$ as the material manifold and model growth by
allowing $\mathcal{B}$ to have an evolving geometry. Consider a
two-dimensional problem, where the relaxed state of the material
is described by a surface. If the bulk of the material grows as,
for example, in a thin shell of biological material undergoing
cell division, the shape of the surface describing the relaxed
state will change. The stresses for a given configuration should
be calculated in terms of the map from the surface describing the
relaxed state, to the current configuration.

While the multiplicative decomposition of the deformation gradient
has been a source of useful approaches to nonlinear problems, we
believe that in many cases, such an approach obfuscates the
underlying natural geometry. A multiplicative decomposition seems
natural if one starts with a stress-free material body and
considers processes such as plasticity which, in general, induce
stresses. However, an initial stress-free Euclidean configuration
may not even exist in certain problems. Mathematically, one can
still consider an incompatible local deformation that brings the
material to a relaxed, Euclidean state, and measure deformations
from this state, as in the multiplicative decomposition described
above. However, we believe a more natural way of looking at the
problem involves treating the material manifold as a non-Euclidean
manifold, and giving its geometry explicitly in terms of the
physics of this problem. In passing we should mention Miehe's
\citep{Miehe1998} work in which instead of
$\mathbf{F}=\mathbf{F}_e\mathbf{F}_p$, he introduces a ``plastic
metric", although with no clear physical meaning/interpretation
for this metric. We will come back to a geometric interpretation
of $\mathbf{F}=\mathbf{F}_e\mathbf{F}_g$ in \S 3.

In this paper, we model growth by introducing a Riemannian
material manifold with an evolving metric. As will be seen, this
formalism is very similar to the approach of
\citep{OzakinYavari2009} to thermoelasticity, however, there are a
few important modifications. First, although we had a version of
mass conservation in thermoelasticity in terms of the changing
material manifold, for the case of growth, mass will in general be
added to (or removed from) the material body; we will have a mass balance. Thus, one can
represent the amount of mass being added (or removed) in terms of
the changes of the differential form describing material
mass-density. Secondly, for the case of thermal stresses, the
material metric was explicitly given in terms of the temperature,
but no such simple dependence exists for the material metric in
biological growth. We will begin by exploring the consequences of
various simple modes of growth, such as a cylindrically symmetric
growth represented by a radius-dependent conformal scaling of the
metric. Assuming simple constitutive relations, we will write the
equations for equilibrium configurations in terms of the
time-dependent metric, much like the case in thermoelasticity. We
will also establish the connection to the formulations involving
multiplicative decomposition of the deformation gradient. 

One should note that there is in fact no guarantee that a
time-dependent Riemannian metric and its Levi-Civita connection is
capable of modeling all kinds of growth. We believe that at the
very least one needs to consider time-dependent connections with
torsion, however, ``Riemannian growth'' is a good starting point.
We aim to investigate the case of growth with torsion, as well as
non-metricity in future communications.

\citet{Efrati2009} have recently studied similar problems in the framework of linearized elasticity by modifying the definition of linearized strain. Here, we start with nonlinear elasticity and instead of modifying any definition of strain will work with an evolving material manifold. We should mention that the idea of using differential geometry in elasticity goes back to more than fifty years ago in the work of \citet{Eckart1948} who realized that the stress-free configuration of a material body evolves in time and an Euclidean stress-free configuration is not always possible. Later developments are due to Kondo \citep{Kondo1955a,Kondo1955b} and Bilby \citep{BilbyGardnerStroh1956,BilbyBulloughSmith1955}.

There have been growth models in the literature using mixture theories. For growth mechanics purposes, a mixture theory is certainly more realistic than a mono-phasic continuum theory. However, in this paper for the sake of simplicity and clarity of presentation, we restrict ourselves to mono-phasic continua. We should mention that our ideas are similar, in spirit, to those of \citet{RajagopalSrinivasa2004b} who have been advocating the idea of material bodies with evolving natural configurations. Here, we work in a fully geometric framework and model a growing body by a continuum that has an evolving Riemannian material metric. One should note that this is a very special case of a possible evolving material manifold that we believe is sufficient for bulk growth purposes. In particular, we work with Levi-Civita connections that are torsion-free. We should also emphasize that we are not, by any means, questioning the usefulness of the traditionally used $\mathbf{F}=\mathbf{F}_e\mathbf{F}_g$ decomposition of deformation gradient. However, we believe that although the existing models based on this multiplicative decomposition have been very useful in growth mechanics (see \citet{BenAmarGoriely2005,Garikipati2004} for some concrete examples.) they all lack a rigorous mathematical foundation. 

This paper is organized as follows. In \S2 we modify the existing
geometric theory of elasticity for growing solids. We show that
energy balance has to be modified and then study its covariance.
We study the entropy production inequality and the restrictions it
imposes on constitutive equations. We also show how Principle of Maximum Entropy Production can be used to obtain thermodynamically-consistent evolution equations for the material metric. We then look at isotropic
growth and model it by a time-dependent rescaling of an initial
material metric. We solve three examples of isotropic and
non-isotropic growth analytically. We then discuss how an evolving material manifold can be visualized using embeddings. In the last part of this
section we obtain stress-free isotropic growth distributions. We then study growth in the Lagrangian field theory of elasticity and show how all the governing equations can be obtained using the Lagrange-d'Alembert principle and using Rayleigh's dissipation functions. In
\S3, we make a connection between the exiting theories of growth
based on the decomposition $\mathbf{F}=\mathbf{F}_e\mathbf{F}_g$
and the present geometric theory. The nonlinear geometric theory
is linearized in \S4. In particular, we obtain those isotropic
growth distributions that are stress free in the linearized
setting. Conclusions are given in \S5.

\section{Evolving Material Metrics and Bulk Growth}

There have been previous works on continuum mechanics formulation
of bodies with variable mass (see \citet{LubardaHoger2002,BenAmarGoriely2005,DiCarloQuiligotti2002,EpsteinMaugin2000}, and references therein). In these works
it is assumed that the growth part of deformation gradient is an unknown tensor field and its evolution is
given by a kinetic equation. In writing energy balance, the
corresponding thermodynamic forces show up. In the present
geometric theory, we work with an evolving material manifold
instead of introducing new fields other than material mass density. Before going into the details
of the proposed theory, let us first briefly review the geometric
theory of classical elasticity.

\paragraph{Geometric Elasticity.} A body $\mathcal{B}$ is identified
with a Riemannian manifold $\mathcal{B}$ and a configuration of
$\mathcal{B}$ is a mapping $\varphi:\mathcal{B} \rightarrow
\mathcal{S}$, where $\mathcal{S}$ is another Riemannian manifold.
The set of all configurations of $\mathcal{B}$ is denoted by
$\mathcal{C}$. A motion is a curve $c:\mathbb{R}\rightarrow
\mathcal{C}; t\mapsto \varphi_t $ in $\mathcal{C}$. It is assumed
that the body is stress free in the material manifold. For a fixed
$t$, $\varphi_t(\mathbf{X})=\varphi(\mathbf{X},t)$ and for a fixed
$\mathbf{X}$, $\varphi_{\mathbf{X}}(t)=\varphi(\mathbf{X},t)$,
where $\mathbf{X}$ is position of material points in the
undeformed configuration $\mathcal{B}$. The material velocity is
the map $\mathbf{V}_t:\mathcal{B}\rightarrow \mathbb{R}^3$ given
by
\begin{equation}
    \mathbf{V}_t(\mathbf{X})=\mathbf{V}(\mathbf{X},t)=\frac{\partial \varphi(\mathbf{X},t)}{\partial
    t}=\frac{d}{dt}\varphi_{\mathbf{X}}(t).
\end{equation}
The material acceleration is defined by
\begin{equation}
    \mathbf{A}_t(\mathbf{X})=\mathbf{A}(\mathbf{X},t)=\frac{\partial \mathbf{V}(\mathbf{X},t)}{\partial
    t}=\frac{d}{dt}\mathbf{V}_{\mathbf{X}}(t).
\end{equation}
In components
\begin{equation}
    A^a=\frac{\partial V^a}{\partial t}+\gamma^a_{bc}V^bV^c,
\end{equation}
where $\gamma^a_{bc}$ is the Christoffel symbol of the local
coordinate chart $\{x^a\}$. Note that $\mathbf{A}$ does not depend
on the connection coefficients of the material manifold.
$\varphi_t$ is assumed to be invertible and regular. The spatial
velocity of a regular motion $\varphi_t$ is defined as
\begin{equation}
    \mathbf{v}_t:\varphi_t(\mathcal{B})\rightarrow
    \mathbb{R}^3,~~~~\mathbf{v}_t=\mathbf{V}_t\circ
    \varphi_t^{-1},
\end{equation}
and the spatial acceleration $\mathbf{a}_t$ is defined as
\begin{equation}
    \mathbf{a}=\dot{\mathbf{v}}=\frac{\partial \mathbf{v}}{\partial
    t}+\mathbf{\nabla}_{\mathbf{v}}\mathbf{v}.
\end{equation}
In components
\begin{equation}
    a^a=\frac{\partial v^a}{\partial t}+\frac{\partial v^a}{\partial x^b}v^b+\gamma^a_{bc}v^bv^c.
\end{equation}

\begin{figure}[hbt]
\begin{center}
\includegraphics[scale=0.65,angle=0]{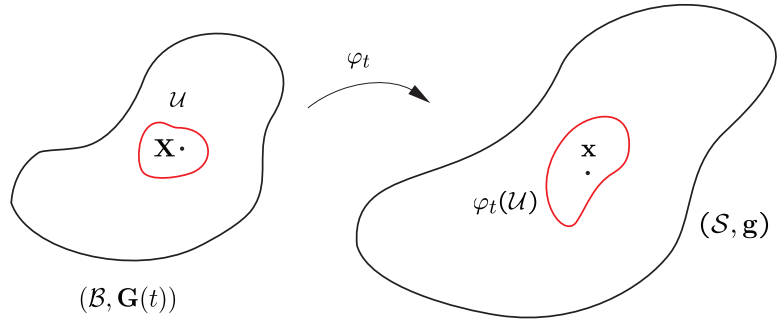}
\end{center}{\vskip -0.1 in}
\caption{\footnotesize Deformation of a growing body as a map between a Riemannian material manifold with a dynamic metric and an ambient space with a fixed background metric.} \label{Deformation}
\end{figure}

Let $\varphi:\mathcal{B}\rightarrow \mathcal{S}$ be a $C^1$
configuration of $\mathcal{B}$ in $\mathcal{S}$, where
$\mathcal{B}$ and $\mathcal{S}$ are manifolds. Deformation gradient is the tangent map of $\varphi$ and is
denoted by $\mathbf{F}=T\varphi$. Thus, at each point $\mathbf{X}
\in \mathcal{B}$, it is a linear map
\begin{equation}
   \mathbf{F}(\mathbf{X}):T_{\mathbf{X}}\mathcal{B}\rightarrow
    T_{\varphi(\mathbf{X})}\mathcal{S}.
\end{equation}
If $\{x^a\}$ and $\{X^A\}$ are local coordinate charts on
$\mathcal{S}$ and $\mathcal{B}$, respectively, the components of
$\mathbf{F} $ are
\begin{equation}
    F^a{}_{A}(\mathbf{X})=\frac{\partial \varphi^a}{\partial X^A}(\mathbf{X}).
\end{equation}
Suppose $\mathcal{B}$ and $\mathcal{S}$ are Riemannian manifolds
with inner products $\left\langle \! \left\langle, \right\rangle
\! \right\rangle_{\mathbf{G}}$ and $\left\langle \! \left\langle,
\right\rangle \! \right\rangle_{\mathbf{g}}$ based at
${\mathbf{X}} \in \mathcal{B}$ and ${\mathbf{x}} \in \mathcal{S}$,
respectively. Transpose of  $\mathbf{F}$ is
defined by
\begin{equation}
    \mathbf{F}^{\textsf{T}}:T_{\mathbf{x}}\mathcal{S}  \rightarrow T_{\mathbf{X}}\mathcal{B},~~~
    \left\langle \! \left\langle \mathbf{FV},\mathbf{v}\right\rangle \! \right\rangle_{\mathbf{g}}= \left\langle \! \left\langle
    \mathbf{V},\mathbf{F}^{\textsf{T}}\mathbf{v}\right\rangle \!
    \right\rangle_{\mathbf{G}}~~~~~\forall~\mathbf{V} \in
    T_{\mathbf{X}}\mathcal{B},~\mathbf{v} \in T_{\mathbf{x}}
    \mathcal{S}.
    \end{equation}
In components
\begin{equation}
    (F^{\textsf{T}}(\mathbf{X}))^A{}_{a}=g_{ab}(\mathbf{x})F^b{}_{B}(\mathbf{X})G^{AB}(\mathbf{X}).
\end{equation}
The right Cauchy-Green deformation tensor is defined by
\begin{equation}
    \mathbf{C}(X):T_{\mathbf{X}} \mathcal{B}\rightarrow T_{\mathbf{X}} \mathcal{B},~~~~~\mathbf{C}(\mathbf{X})=\mathbf{F}(\mathbf{X})^{\textsf{T}} \mathbf{F}(\mathbf{X}),
\end{equation}
where $\mathbf{g}$ and $\mathbf{G}$ are metric tensors on
$\mathcal{S}$ and $\mathcal{B}$, respectively. In components
\begin{equation}
    C^A_{~B}=(F^{\textsf{T}})^A{}_{a}F^a{}_{B}.
\end{equation}
One can show that
\begin{equation}
    \mathbf{C}^\flat=\varphi^*(\mathbf{g})=\mathbf{F}^*\mathbf{g}\mathbf{F},~\textrm{i.e.}~~~C_{AB}=(g_{ab}\circ
    \varphi)F^a{}_{A}F^b{}_{B}.
\end{equation}
For bulk growth, we assume that the material manifold
$\mathcal{B}$ remains unchanged but the metric evolves, i.e.
$\mathbf{G}=\mathbf{G}(\mathbf{X},t)$\footnote{In mathematics,
evolving metrics have been studied extensively. The most
celebrated example is Ricci flow \citep{Hamilton1982,Topping2006} that was
used in proving Poincar\'e Conjecture by \citet{Perelman2002}.
Interestingly, for a seemingly very different application, i.e.
growth mechanics, an evolving geometry plays a key role.}. When
mass is added or removed, the stress-free state of the body
changes. Local changes in mass change the stress-free configuration of the body. This is modeled by a
time-dependent material metric that represents local changes in
volume and ``shape" in the relaxed configuration (see Fig. \ref{Deformation}). In \S3, we will
make a connection between this approach and the conventional
$\mathbf{F}=\mathbf{F}_e\mathbf{F}_p$ decomposition of deformation
gradient.

\paragraph{Incompressibility.} In growth mechanics it is usually assumed that elastic deformations are incompressible. In the classical theory in which $\mathbf{F}=\mathbf{F}_e\mathbf{F}_g$ is assumed, incompressibility implies $\det \mathbf{F}_e =1$, i.e. all the volume changes are due to growth. In the geometric theory the following relation holds between volume elements of $(\mathcal{B},\mathbf{G})$ and $(\mathcal{S},\mathbf{g})$:
\begin{equation}
    dv= J~dV,
\end{equation}
where
\begin{equation}
    J=\sqrt{\frac{\det \mathbf{g}}{\det \mathbf{G}}}\det \mathbf{F}.
\end{equation}
Incompressibility of elastic deformations means that $J=1$. Note that even when $J=1$, still $dv$ is time dependent as a result of the time evolution of the material metric that makes $dV$ time dependent. In other words, an observer in the ambient space sees changes in volume that are only due to volume changes in the material manifold. We will show the equivalence of $J_e=1$ in the classical theory with $J=1$ in the geometric theory in both some simple examples in \S2.9 and in the general case in \S3.

\subsection{Energy Balance.} Let us look at energy balance for a growing body. The standard
material balance of energy for a subset
$\mathcal{U}\subset\mathcal{B}$ reads \citep{YaMaOr2006}
\begin{equation}\label{Energy-Balance-Standard}
    \frac{d}{d
    t}\int_{\mathcal{U}}\rho_0\left(E+\frac{1}{2}\left\langle \! \left\langle\mathbf{V},\mathbf{V}\right\rangle \! \right\rangle\right)dV=
    \int_{\mathcal{U}}\rho_0\left(\left\langle \! \left\langle\mathbf{B},\mathbf{V}\right\rangle \! \right\rangle+R\right)dV+
    \int_{\partial\mathcal{U}}\left(\left\langle \! \left\langle\mathbf{T},\mathbf{V}\right\rangle \!
    \right\rangle+H\right)dA,
\end{equation}
where
$E=E(\mathbf{X},\textsf{N},\mathbf{G},\mathbf{F},\mathbf{g}\circ\varphi)$ is
the material internal energy density, $\textsf{N}$, $\rho_0$, $\mathbf{B}$,
$\mathbf{T}$, $R$, and $H$ are specific entropy, material mass density, body force
per unit undeformed mass, traction vector, heat supply, and heat
flux, respectively.

We first note that energy balance should be modified in the case
of growing bodies with time-dependent material metrics. Note that
when metric is time dependent, material density mass form
$\textsf{m}(\mathbf{X},t)=\rho_0(\mathbf{X},t)dV(\mathbf{X},t)$ is
time dependent even if $\rho_0$ is not time dependent. For a
subbody $\mathcal{U}\subset\mathcal{B}$, rate of change of mass
reads
\begin{equation}
    \frac{d}{d t}\int_{\mathcal{U}}\rho_0(\mathbf{X},t)dV(\mathbf{X},t)
    =\int_{\mathcal{U}}\left[\frac{\partial \rho_0}{\partial t}+\frac{1}{2}\rho_0\operatorname{tr}\!\left(\frac{\partial \mathbf{G}}{\partial t}\right)\right]dV.
\end{equation}
Note that if $\rho_0$ is time independent, then the term
$\frac{1}{2}\rho_0\operatorname{tr}\!\left(\frac{\partial
\mathbf{G}}{\partial t}\right)$ represents the change in mass due
to growth. Here, we assume that a scalar field of mass source/sink $S_m(\mathbf{X},t)$ is given.\footnote{Note that by definition
\begin{equation*}
    \frac{d}{dt}\int_{\mathcal{U}}\rho_0(X,t)dV=\int_{\mathcal{U}}S_m\!(X,t) dV=\int_{\mathcal{U}}\stackrel{\tiny \circ}{S}_m(X,t) d\!\stackrel{\tiny \circ}{V},
\end{equation*}
where $\stackrel{\tiny \circ}{S}_m\!(X,t)$ is mass source in the initial material manifold with volume element $d\!\stackrel{\tiny \circ}{V}$. Note also that physically $\stackrel{\tiny \circ}{S}$ is given.} This mass source will change the stress-free configuration of the body and $(\mathcal{B},\mathbf{G}(\mathbf{X},t))$ represents the stress-free configuration of the body.

The rate of change of material metric is a kinematical variable
that contributes to power. Therefore, energy balance for a growing
body with a time-dependent material metric is written
as\footnote{Note that in \citet{LubardaHoger2002} the term
analogous to $\frac{\partial \rho_0}{\partial
t}+\frac{1}{2}\rho_0\operatorname{tr}\!\left(\frac{\partial
\mathbf{G}}{\partial t}\right)$ is denoted by $r_g$. There,
instead of the term $\rho_0\frac{\partial E}{\partial \mathbf{G}}:\frac{\partial\mathbf{G}}{\partial t}$ they introduce
a term $\rho \mathcal{R}_g r_g$. One should note that even if mass
is conserved at a point, still a change in shape can contribute to
energy balance and is captured in our formulation. See also \citet{EpsteinMaugin2000} and \citet{LubardaHoger2002}.}
\begin{eqnarray}
  && \frac{d}{d
    t}\int_{\mathcal{U}}\rho_0\left(E+\frac{1}{2}\left\langle \! \left\langle\mathbf{V},\mathbf{V}\right\rangle \! \right\rangle\right)dV
    =\int_{\mathcal{U}}\Bigg\{\rho_0\left(\left\langle \! \left\langle\mathbf{B},\mathbf{V}\right\rangle \! \right\rangle+R\right)
    +\rho_0\frac{\partial E}{\partial \mathbf{G}}:\frac{\partial\mathbf{G}}{\partial t}
    \nonumber\\
  \label{Energy-Balance-Growth} && ~~~~~~~~~~~~~~~~~~~~~~~~~~~~~~~~~~~~~~~~~~~~+S_m  \left(\! E+\frac{1}{2}\left\langle \! \left\langle\mathbf{V},\mathbf{V}\right\rangle \! \right\rangle \!\right)
  \Bigg\}
    dV +\int_{\partial\mathcal{U}}\left(\left\langle \! \left\langle\mathbf{T},\mathbf{V}\right\rangle \!
    \right\rangle+H\right)dA.
\end{eqnarray}

\subsection{Covariance of Energy Balance.} It turns out that in continuum mechanics (and even discrete systems) one can
obtain all the balance laws using energy balance and postulating
its invariance under some groups of transformations. This was
introduced by \citet{Green64} in the case of Euclidean ambient
spaces and was extended to manifolds by \citet{MaHu1983}. See also
\citet{SiMa1984,YaMaOr2006,YavariOzakin2008,Yavari2008,YavariMarsden2009a,YavariMarsden2009b}
for applications of covariance ideas in different continuous and
discrete systems.

In order to covariantly obtain all the balance laws, we postulate
that energy balance is form invariant under an arbitrary
time-dependent spatial diffeomorphism
$\xi_t:\mathcal{S}\rightarrow\mathcal{S}$, i.e.
\begin{eqnarray}
  && \frac{d}{d
    t}\int_{\mathcal{U}}\rho'_0\left(E'+\frac{1}{2}\left\langle \! \left\langle\mathbf{V}',\mathbf{V}'\right\rangle \! \right\rangle\right)dV
    =\int_{\mathcal{U}}\Bigg\{\rho'_0\left(\left\langle \! \left\langle\mathbf{B}',\mathbf{V}'\right\rangle \! \right\rangle+R'\right)
    +\rho'_0\frac{\partial E'}{\partial \mathbf{G}'}:\frac{\partial\mathbf{G}'}{\partial t}
     \nonumber\\
  \label{Energy-Balance-Growth-primed} && ~~~~~~~~~~~~~~~~~~~~~~~~~~~~~~~~~~~~~~~~~~~~+S_m'
  \left(\!E'+\frac{1}{2}\left\langle \! \left\langle\mathbf{V}',\mathbf{V}'\right\rangle \! \right\rangle\!\right)
  \Bigg\}
    dV +\int_{\partial\mathcal{U}}\left(\left\langle \! \left\langle\mathbf{T}',\mathbf{V}'\right\rangle \!
    \right\rangle+H'\right)dA.
\end{eqnarray}
Note that \citep{YaMaOr2006}
\begin{equation}
    R'=R,~~H'=H,~~\rho'_0=\rho_0,~~\mathbf{T}'=\xi_{t*}\mathbf{T},~~\mathbf{V}'=\xi_{t*}\mathbf{V}+\mathbf{W},
\end{equation}
where $\mathbf{W}=\frac{\partial}{\partial t}\xi_t\circ\varphi$.
Note also that
\begin{equation}
    \mathbf{G}'=\mathbf{G},\frac{\partial \mathbf{G}'}{\partial
  t}=\frac{\partial \mathbf{G}}{\partial t}~~~~~\textrm{and}~~~~~
  E'(\mathbf{X},\textsf{N}',\mathbf{G},\mathbf{F}',\mathbf{g}\circ\varphi')=E(\mathbf{X},\textsf{N},\mathbf{G},\mathbf{F},\xi_t^*\mathbf{g}\circ\varphi).
\end{equation}
Thus, at $t=t_0$
\begin{equation}
    \frac{d}{dt}E'=\frac{\partial E}{\partial \textsf{N}}:\frac{d\textsf{N}}{d t}
    +\frac{\partial E}{\partial\mathbf{G}}:\frac{\partial\mathbf{G}}{\partial t}+\frac{\partial E}{\partial\mathbf{g}\circ\varphi}:\mathfrak{L}_{W}\mathbf{g}\circ\varphi.
\end{equation}
We also assume that body forces are transformed such that
\citep{MaHu1983}
$\mathbf{B}'-\mathbf{A}'=\xi_{t*}(\mathbf{B}-\mathbf{A})$.
Therefore, (\ref{Energy-Balance-Growth-primed}) at $t=t_0$ reads
\begin{eqnarray}
  && \int_{\mathcal{U}}\left[\frac{\partial \rho_0}{\partial t}+\frac{1}{2}\rho_0\operatorname{tr}\!\left(\frac{\partial \mathbf{G}}{\partial t}\right)\right]
    \left(E+\frac{1}{2}\left\langle \! \left\langle\mathbf{V}+\mathbf{W},\mathbf{V}+\mathbf{W}\right\rangle \!
    \right\rangle\right)dV \nonumber\\
  && ~+\int_{\mathcal{U}} \rho_0\left(\frac{\partial E}{\partial\mathbf{G}}:\frac{\partial\mathbf{G}}{\partial t}+\frac{\partial E}{\partial\mathbf{g}\circ\varphi}:\mathfrak{L}_{W}\mathbf{g}\circ\varphi
    + \left\langle \! \left\langle\mathbf{V}+\mathbf{W},\mathbf{A}\right\rangle \! \right\rangle \right) dV
     \nonumber \\
    &&  =\int_{\mathcal{U}}\left\{\rho_0\left(\left\langle \! \left\langle\mathbf{B},\mathbf{V}+\mathbf{W}\right\rangle \! \right\rangle+R\right)
    +\rho_0\frac{\partial E}{\partial \mathbf{G}}:\frac{\partial\mathbf{G}}{\partial t}
    +S_m\left(E+\frac{1}{2}\left\langle \! \left\langle\mathbf{V}+\mathbf{W},\mathbf{V}+\mathbf{W}\right\rangle \! \right\rangle\right)
    \right\}dV \nonumber\\
   \label{Energy-Balance-Growth-primed-simplified}&& ~~~ +\int_{\partial\mathcal{U}}\left(\left\langle \! \left\langle\mathbf{T},\mathbf{V}+\mathbf{W}\right\rangle \!\right\rangle+H\right)dA.
\end{eqnarray}
Subtracting (\ref{Energy-Balance-Growth}) from
(\ref{Energy-Balance-Growth-primed-simplified}), one obtains
\begin{eqnarray}
  && \int_{\mathcal{U}} \left[\frac{\partial \rho_0}{\partial t}+\frac{1}{2}\rho_0\operatorname{tr}\!\left(\frac{\partial \mathbf{G}}{\partial t}\right)-S_m\right]
    \left(\frac{1}{2}\left\langle \! \left\langle \mathbf{W},\mathbf{W}\right\rangle \!
    \right\rangle +\left\langle \! \left\langle\mathbf{V},\mathbf{W}\right\rangle \!
    \right\rangle\right)  dV    \nonumber \\
   &&  ~~~+
    \int_{\mathcal{U}} \rho_0\left(\frac{\partial E}{\partial\mathbf{g}\circ\varphi}:\mathfrak{L}_{W}\mathbf{g}\circ\varphi
    + \left\langle \! \left\langle\mathbf{A},\mathbf{W}\right\rangle \! \right\rangle \right) dV
    =\int_{\mathcal{U}}\rho_0\left(\left\langle \! \left\langle\mathbf{B},\mathbf{W}\right\rangle \! \right\rangle\right)
    dV +\int_{\partial\mathcal{U}}\left\langle \! \left\langle\mathbf{T},\mathbf{W}\right\rangle \!\right\rangle dA.
\end{eqnarray}
From this and arbitrariness of $\mathbf{W}$ we conclude that
\citep{YaMaOr2006}
\begin{eqnarray}
  && \frac{\partial \rho_0}{\partial t}+\frac{1}{2}\rho_0\operatorname{tr}\!\left(\frac{\partial \mathbf{G}}{\partial t}\right)=S_m,\\
  && \operatorname{Div}\mathbf{P}+\rho_0\mathbf{B}=\rho_0\mathbf{A},  \\
  && \label{DE-formula} 2\rho_0\frac{\partial E}{\partial\mathbf{g}\circ\varphi}= \boldsymbol{\tau},  \\
  && \boldsymbol{\tau}^{\textsf{T}}=\boldsymbol{\tau},
\end{eqnarray}
where $\mathbf{P}$ is the first Piola-Kirchhoff stress and $\boldsymbol{\tau}=J\boldsymbol{\sigma}$ is the Kirchhoff
stress. It is seen that instead of conservation of mass we have a balance of mass and the remaining
balance laws are unchanged. Note, however, that divergence and
acceleration both explicitly depend on $\mathbf{G}$, i.e. the time
dependency of material metric affects the governing balance
equations. We will see examples in \S2.9.

\subsection{Local Form of Energy Balance.} Let us now localize the
energy balance. First note that
\begin{equation}
    \frac{d}{d t}E=\mathbf{L}_{V}E=\frac{\partial E}{\partial \textsf{N}}\frac{d \textsf{N}}{d t}+\frac{\partial E}{\partial \mathbf{G}}:\frac{\partial\mathbf{G}}{\partial t}
    +\frac{\partial E}{\partial \mathbf{F}}:\mathbf{L}_{\mathbf{V}}\mathbf{F}
    +\frac{\partial E}{\partial \mathbf{g}}:\mathbf{L}_{\mathbf{V}}\mathbf{g}\circ\varphi.
\end{equation}
Note that $\mathbf{L}_{\mathbf{V}}\mathbf{F}=\mathbf{0}$ because
for an arbitrary $\mathbf{Z}\in T_{\mathbf{X}}\mathcal{B}$
\begin{equation}
    \mathbf{L}_{\mathbf{V}}\mathbf{F}=\frac{\partial}{\partial  t}\varphi^*\left(\mathbf{F}\cdot\mathbf{Z}\right)
    =\frac{\partial}{\partial  t}\varphi^*\left(\varphi_*\mathbf{Z}\right)
    =\frac{\partial}{\partial  t}\mathbf{Z}=\mathbf{0}.
\end{equation}
Using this and also noting that because the background metric is
time independent, we have
\begin{equation}
    \frac{d}{d t}E=\frac{\partial E}{\partial \textsf{N}}\frac{d \textsf{N}}{d t}+ \frac{\partial E}{\partial \mathbf{G}}:\frac{\partial\mathbf{G}}{\partial t}
    +\frac{\partial E}{\partial \mathbf{g}}:\mathbf{d},
\end{equation}
where
$\mathbf{d}=\frac{1}{2}\mathfrak{L}_{\mathbf{V}}\mathbf{g}\circ\varphi$
is the rate of deformation tensor.\footnote{This is the symmetric part of $\nabla \mathbf{v}$, i.e. the symmetric part of the so-called ``velocity gradient".} We know that $H=-\left\langle
\!\! \left\langle\mathbf{Q},\hat{\mathbf{N}}\right\rangle \!\!
\right\rangle$ and \citep{YaMaOr2006}
\begin{equation}
    \int_{\partial\mathcal{U}}\left\langle \! \left\langle\mathbf{T},\mathbf{V}\right\rangle \!
    \right\rangle dA=\int_{\mathcal{U}}\left(\left\langle \! \left\langle\operatorname{Div}\mathbf{P},\mathbf{V}\right\rangle \!
    \right\rangle +\boldsymbol{\tau}:\boldsymbol{\Omega} + \boldsymbol{\tau}:\mathbf{d}
    \right)dV,
\end{equation}
where $\Omega_{ab}=\frac{1}{2}(V_{a|b}-V_{b|a})$, and $\boldsymbol{\tau}$ is Kirchhoff stress. Thus, from
(\ref{Energy-Balance-Growth}) and using balances of linear and angular momenta we obtain the local form of energy
balance as
\begin{equation}
    \rho_0\frac{d E}{dt}+\operatorname{Div}\mathbf{Q}=\rho_0\frac{\partial E}{\partial \mathbf{G}}:\frac{\partial\mathbf{G}}{\partial t}
    +\boldsymbol{\tau}:\mathbf{d}+\rho_0 R.
\end{equation}
In term of the first Piola-Kirchhoff stress, this can be written
as
\begin{equation}\label{EnergyBalance-Local}
    \rho_0\frac{d E}{dt}+\operatorname{Div}\mathbf{Q}=\rho_0\frac{\partial E}{\partial \mathbf{G}}:\frac{\partial\mathbf{G}}{\partial t}
    +\mathbf{P}:\nabla_0\mathbf{V}+\rho_0 R,
\end{equation}
where
$\mathbf{P}:\nabla_0\mathbf{V}=P^{aA}V^a{}_{|A}$.

\paragraph{Material Metric Evolution.} Evolution of material metric is assumed to be given through a kinetic
equation of the form\footnote{A simple example that has been
extensively studied is the Ricci flow \citep{Hamilton1982}, which
reads
\begin{equation*}
    \frac{\partial\mathbf{G}}{\partial  t}=-2\mathbf{R},
\end{equation*}
where $\mathbf{R}$ is the Ricci curvature of the metric
$\mathbf{G}$. Note that this flow smooth outs a distorted 3-sphere
\citep{Topping2006}.}
\begin{equation}
    \frac{\partial\mathbf{G}}{\partial
    t}=\boldsymbol{\Upsilon}(\mathbf{X},\mathbf{G},\mathbf{F},\mathbf{g})
    =\boldsymbol{\Phi}(\mathbf{X},\mathbf{G},\mathbf{C}).
\end{equation}
See \citet{AmbrosiMollica2004,LoretSimoes2005,Fusi2006,AmbrosiGuana2007} for some examples written in terms
of evolution of $\mathbf{F}_g$. We will come back to this problem after first discussing the Second Law of Thermodynamics for a growing body.

\subsection{The Second Law of Thermodynamics and Restrictions on Constitutive Equations}

In the absence of growth, entropy production inequality in
material coordinates has the following form
\citep{ColemanNoll1963}
\begin{equation}
    \frac{d}{dt}\int_{\mathcal{U}}\rho_0 \textsf{N}dV \geq  \int_{\mathcal{U}}\frac{\rho_0
    R}{\Theta}dV+\int_{\partial\mathcal{U}}\frac{H}{\Theta}~dA,
\end{equation}
where
$\textsf{N}=\textsf{N}(\mathbf{X},t)$
is the material entropy density and $\Theta=\Theta(\mathbf{X},t)$
is the absolute temperature. This is called the Clausius-Duhem
inequality.

When material metric is time dependent using balance of mass the Clausius-Duhem
inequality is modified to read\footnote{See also \citet{EpsteinMaugin2000}.}
\begin{equation}\label{Entropy}
    \frac{d}{dt}\int_{\mathcal{U}}\rho_0 \textsf{N}dV \geq  \int_{\mathcal{U}}\frac{\rho_0
    R}{\Theta}dV+\int_{\partial\mathcal{U}}\frac{H}{\Theta}~dA
    +\int_{\mathcal{U}} \textsf{N} \left[\frac{\partial \rho_0}{\partial t}+\frac{1}{2}\rho_0\operatorname{tr}\!\left(\frac{\partial \mathbf{G}}{\partial
    t}\right)\right] dV+\int_{\mathcal{U}} \rho_0\frac{\partial E}{\partial \mathbf{G}}:\frac{\partial\mathbf{G}}{\partial t}
    dV.
\end{equation}
This inequality can be localized to read
\begin{equation}\label{Entropy-Local}
    \rho_0 \frac{d \textsf{N}}{dt} \geq  \frac{\rho_0
    R}{\Theta}-\operatorname{Div}\!\left(\frac{\mathbf{Q}}{\Theta}\right)
    + \rho_0\frac{\partial E}{\partial \mathbf{G}}:\frac{\partial\mathbf{G}}{\partial t}.
\end{equation}
Note that $E=\Psi+\textsf{N}\Theta$ and hence\footnote{Note that $\Psi=\Psi(\mathbf{X},\Theta,\mathbf{G},\mathbf{F},\mathbf{g}\circ\varphi)$ and $E=E(\mathbf{X},\textsf{N},\mathbf{G},\mathbf{F},\mathbf{g}\circ\varphi)$.}
\begin{equation}
    \frac{d}{dt}\textsf{N}=\frac{1}{\Theta}\left(\frac{d E}{dt}-\frac{d \Psi}{dt}\right)-\frac{\dot{\Theta}}{\Theta^2}(E-\Psi).
\end{equation}
Substituting this into (\ref{Entropy-Local}) yields
\begin{equation}
    \frac{\rho_0}{\Theta}\frac{d E}{dt}-\frac{\rho_0}{\Theta}\frac{d \Psi}{dt}-\rho_0\frac{\dot{\Theta}}{\Theta^2}(E-\Psi) \geq  \frac{\rho_0
    R}{\Theta}-\operatorname{Div}\!\left(\frac{\mathbf{Q}}{\Theta}\right)
    + \rho_0\frac{\partial E}{\partial \mathbf{G}}:\frac{\partial\mathbf{G}}{\partial t}.
\end{equation}
Now substituting the local energy balance
(\ref{EnergyBalance-Local}) into the above inequality we obtain
\begin{equation}\label{CD-Inequality}
    \mathbf{P}:\nabla_0\mathbf{V}-\rho_0\frac{d \Psi}{dt}-\rho_0\textsf{N}\dot{\Theta}
    \geq  \frac{1}{\Theta}\mathbf{d}\Theta\cdot\mathbf{Q}.
\end{equation}
We know that $\Psi=\Psi(\mathbf{X},\Theta,\mathbf{G},\mathbf{F},\mathbf{g}\circ\varphi)$
and thus\footnote{Nota that because $\Psi$ is a scalar its time derivative is
equal to its Lie derivative along the velocity vector field, and
thus when the material metric is time independent we have
\begin{equation*}
    \frac{d}{dt}\Psi=\mathbf{L}_{\mathbf{V}}\Psi=\frac{\partial \Psi}{\partial \Theta}\dot{\Theta}+
    \frac{\partial \Psi}{\partial \mathbf{G}}:\frac{\partial\mathbf{G}}{\partial t}+
    \frac{\partial \Psi}{\partial\mathbf{F}}:\mathbf{L}_{\mathbf{V}}\mathbf{F}+
    \frac{\partial\Psi}{\partial\mathbf{g}}:\mathbf{L}_{\mathbf{V}}\mathbf{g}
    =\frac{\partial \Psi}{\partial \Theta}\dot{\Theta}
    +\frac{\partial \Psi}{\partial \mathbf{G}}:\frac{\partial\mathbf{G}}{\partial t}
    +\frac{\partial\Psi}{\partial\mathbf{g}}:\mathbf{L}_{\mathbf{V}}\mathbf{g}.
\end{equation*}
We also know that the same time derivative is equal to the
covariant derivative of \textsf{N} with respect to velocity
vector,i.e.
\begin{equation*}
    \frac{d}{dt}\Psi=\frac{\partial}{\partial t}\Psi+\nabla_{\mathbf{V}}\Psi=\frac{\partial \Psi}{\partial \Theta}\dot{\Theta}+
    \frac{\partial \Psi}{\partial \mathbf{G}}:\frac{\partial\mathbf{G}}{\partial t}+
    \frac{\partial \Psi}{\partial\mathbf{F}}:\nabla_{\mathbf{V}}\mathbf{F}+
    \frac{\partial\Psi}{\partial\mathbf{g}}:\nabla_{\mathbf{V}}\mathbf{g}
    =\frac{\partial \Psi}{\partial \Theta}\dot{\Theta}+
    \frac{\partial \Psi}{\partial \mathbf{G}}:\frac{\partial\mathbf{G}}{\partial t}
    +\frac{\partial \Psi}{\partial\mathbf{F}}:\nabla_{\frac{\partial}{\partial t}}{\mathbf{F}}
    =\frac{\partial \Psi}{\partial \Theta}\dot{\Theta}
    +\frac{\partial \Psi}{\partial \mathbf{G}}:\frac{\partial\mathbf{G}}{\partial t}
    +\frac{\partial \Psi}{\partial\mathbf{F}}:\nabla_0\mathbf{V},
\end{equation*}
where we used the fact that $\frac{\partial}{\partial t}\Psi=0$ and $\nabla_{\frac{\partial}{\partial
t}}{\mathbf{F}}=\nabla_0\mathbf{V}$. See \citet{Nishikawa2002} for a proof.
}
\begin{equation}
    \frac{d \Psi}{d t}=\frac{\partial \Psi}{\partial\Theta}\dot{\Theta}
    +\frac{\partial \Psi}{\partial\mathbf{G}}:\frac{\partial\mathbf{G}}{\partial t}
    +\frac{\partial\Psi}{\partial \mathbf{F}}:\nabla_0\mathbf{V}.
\end{equation}
Therefore, (\ref{CD-Inequality}) is simplified to read
\begin{equation}
    \rho_0\left(\frac{\partial \Psi}{\partial\Theta}+\textsf{N}\right)\dot{\Theta}
    +\left(\rho_0\frac{\partial \Psi}{\partial\mathbf{F}}-\mathbf{P}\right):\nabla_0\mathbf{V}
    +\frac{1}{\Theta}\mathbf{d}\Theta\cdot\mathbf{Q}
    +\rho_0\frac{\partial \Psi}{\partial\mathbf{G}}:\dot{\mathbf{G}} \leq 0.
\end{equation}
Following \cite{ColemanNoll1963} and \cite{MaHu1983} we conclude
that
\begin{equation}\label{Restrictions}
    \frac{\partial \Psi}{\partial\Theta}=-\textsf{N}~~~~~\textrm{and}~~~~~\rho_0\frac{\partial
    \Psi}{\partial\mathbf{F}}=\mathbf{P},
\end{equation}
and entropy production inequality reduces to
$\rho_0\frac{\partial \Psi}{\partial\mathbf{G}}:\dot{\mathbf{G}}+\frac{1}{\Theta}\mathbf{d}\Theta\cdot\mathbf{Q} \leq 0$.

\paragraph{Remark.} There have been objections in the literature
on using the Clausius-Duhem inequality in continuum mechanics
\citep{GreenNaghdi1977,MaHu1983}. Next, we show that in growth
mechanics energy balance and a more general notion of covariance
are enough to obtain the restrictions (\ref{Restrictions}) on
constitutive equations. In passing we should mention that
\citet{GreenNaghdi1991} were able to obtain their entropy balance
using energy balance and invariance arguments in the case of
Euclidean ambient space. What we will show next is consistent with
their results.

\subsection{Restrictions on Constitutive Equations Using a Thermomechanical Covariance of Energy Balance}

In this subsection, we follow \citet{MaHu1983} and obtain the
restrictions (\ref{Restrictions}) on the constitutive equations
using covariance of energy balance with no reference to the
entropy production inequality. In the case of classical
elasticity, \citet{MaHu1983} started with the local form of energy
balance and postulated its covariance under simultaneous action of
time-dependent spatial diffeomorphisms and time-dependent
monotonically increasing temperature rescalings. Here we start
with the integral form of the energy balance.

Let us consider spatial diffeomorphisms
$\xi_t:\mathcal{S}\rightarrow\mathcal{S}$ and monotonically
increasing temperature rescalings $\zeta_t:\mathbb{R}^+\rightarrow
\mathbb{R}^+$. We assume that at $t=t_0$, $\zeta=1$ and
$\frac{d}{dt}\zeta_t=z$. We also assume that under these
transformations, energy balance is invariant, i.e.
\begin{eqnarray}
  && \frac{d}{d
    t}\int_{\mathcal{U}}\rho'_0\left(E'+\frac{1}{2}\left\langle \! \left\langle\mathbf{V}',\mathbf{V}'\right\rangle \! \right\rangle\right)dV
    =\int_{\mathcal{U}}\Bigg\{\rho'_0\left(\left\langle \! \left\langle\mathbf{B}',\mathbf{V}'\right\rangle \! \right\rangle+R'\right)
    +\rho'_0\frac{\partial E'}{\partial \mathbf{G}}:\frac{\partial\mathbf{G}}{\partial t}
    +S'_m
  \left(\!E'+\frac{1}{2}\left\langle \! \left\langle\mathbf{V}',\mathbf{V}'\right\rangle \! \right\rangle\!\right)
  \Bigg\}
    dV
     \nonumber\\
   \label{Energy-Balance-Growth-primed-thermal} &&
   ~~~~~~~~~~~~~~~~~~~~~~~~~~~~~~~~~~~~~~~~~~~~~
   +\int_{\partial\mathcal{U}}\left(\left\langle \! \left\langle\mathbf{T}',\mathbf{V}'\right\rangle \!
    \right\rangle+H'\right)dA.
\end{eqnarray}
Note that $E=\Psi+\Theta \textsf{N}$. In the new frame
$\varphi'=\xi\circ\varphi$ and $\Theta'=\zeta \Theta$. We assume
that $E$ transforms tonsorially, i.e.
\begin{equation}
    E'(\mathbf{X},\textsf{N}',\mathbf{G},\mathbf{F}',\mathbf{g}')=E(\mathbf{X},\textsf{N},\mathbf{G},\xi_{*}\mathbf{F},\xi_{*}\mathbf{g}).
\end{equation}
The same transformation is assumed for free energy density and
hence
\begin{equation}
    \frac{d}{dt}\Psi'=\frac{d}{dt}\Psi(\mathbf{X},\zeta\Theta,\mathbf{G},\xi_{*}\mathbf{F},\xi_{*}\mathbf{g})
    =\frac{\partial \Psi}{\partial \Theta'}\dot{\Theta}'+\frac{\partial \Psi}{\partial\mathbf{G}}:\dot{\mathbf{G}}
    +\frac{\partial \Psi}{\partial\xi_{*}\mathbf{F}}:\nabla_0\mathbf{V}',
\end{equation}
where $(\nabla_0\mathbf{V})^a{}_A=V^a{}_{|A}$.
Therefore, at $t=t_0$
\begin{equation}
    \frac{d}{dt}\Psi'=\frac{\partial \Psi}{\partial \Theta}(\dot{\Theta}+z\Theta)+\frac{\partial \Psi}{\partial\mathbf{G}}:\dot{\mathbf{G}}
    +\frac{\partial\Psi}{\partial\mathbf{F}}:(\nabla_0\mathbf{V}+\nabla_0\mathbf{W}).
\end{equation}
Thus, at $t=t_0$, we can write
\begin{equation}
    \frac{d}{dt}E'=\frac{d}{dt}E +\frac{\partial \Psi}{\partial\mathbf{F}}:\nabla_0\mathbf{W}
    +z\left(\frac{\partial \Psi}{\partial \Theta}+\textsf{N}\right)\Theta+\left(\frac{d \textsf{N}'}{dt}\Theta'-\frac{d \textsf{N}}{dt}\Theta\right)_{t=t_0}.
\end{equation}
Energy balance in the new frame at $t=t_0$ is simplified to read
\begin{eqnarray}
  && \int_{\mathcal{U}}\left[ \frac{\partial \rho_0}{\partial t} +\frac{1}{2}\rho_0\operatorname{tr}\left(\frac{\partial \mathbf{G}}{\partial t}\right)\right]
  \left(E+\frac{1}{2}\left\langle \! \left\langle\mathbf{V},\mathbf{V}\right\rangle \! \right\rangle
  +\left\langle \! \left\langle\mathbf{V},\mathbf{W}\right\rangle \! \right\rangle
  +\frac{1}{2}\left\langle \! \left\langle\mathbf{W},\mathbf{W}\right\rangle \! \right\rangle\right)
    dV \nonumber\\
  && ~~+\int_{\mathcal{U}}\rho_0\left[ \frac{d}{dt}E +\frac{\partial \Psi}{\partial\mathbf{F}}:\nabla_0\mathbf{W}
    +z\left(\frac{\partial \Psi}{\partial \Theta}+\textsf{N}\right)\Theta+\left(\frac{d \textsf{N}'}{dt}\Theta'-\frac{d \textsf{N}}{dt}\Theta\right)_{t=t_0}
    \!\!+ \left\langle \! \left\langle\mathbf{V}+\mathbf{W},\mathbf{A}'|_{t=t_0}\right\rangle \!  \right\rangle  \right]dV \nonumber\\
  && ~~~ =\int_{\mathcal{U}}\Bigg[\rho_0\left(\left\langle \! \left\langle\mathbf{B}'|_{t=t_0},\mathbf{V}+\mathbf{W}\right\rangle \! \right\rangle+R'|_{t=t_0}\right)
    +\rho_0\frac{\partial E}{\partial \mathbf{G}}:\frac{\partial\mathbf{G}}{\partial t} \nonumber\\
   && ~~~~~ +S_m \left(E+\frac{1}{2}\left\langle \! \left\langle\mathbf{V},\mathbf{V}\right\rangle \! \right\rangle
  +\left\langle \! \left\langle\mathbf{V},\mathbf{W}\right\rangle \! \right\rangle
  +\frac{1}{2}\left\langle \! \left\langle\mathbf{W},\mathbf{W}\right\rangle \! \right\rangle\right)\Bigg]dV
  +\int_{\partial\mathcal{U}}\left(\left\langle \! \left\langle\mathbf{T},\mathbf{V}+\mathbf{W}\right\rangle \!\right\rangle+H'|_{t=t_0}\right)dA.
\end{eqnarray}
Assuming that
$\mathbf{B}'-\mathbf{A}'=\xi_{*}(\mathbf{B}-\mathbf{A})$ and
subtracting balance of energy (\ref{Energy-Balance-Growth}) from
(\ref{Energy-Balance-Growth-primed-thermal}), we obtain
\begin{eqnarray}
  && \int_{\mathcal{U}}\left[ \frac{\partial \rho_0}{\partial t} +\frac{1}{2}\rho_0\operatorname{tr}\left(\frac{\partial \mathbf{G}}{\partial t}\right)-S_m\right]
  \left(\left\langle \! \left\langle\mathbf{V},\mathbf{W}\right\rangle \! \right\rangle
  +\frac{1}{2}\left\langle \! \left\langle\mathbf{W},\mathbf{W}\right\rangle \! \right\rangle\right)
    dV \nonumber\\
  && ~+\int_{\mathcal{U}}\rho_0\left[ \frac{\partial \Psi}{\partial\mathbf{F}}:\nabla_0\mathbf{W}
    +z\left(\frac{\partial \Psi}{\partial \Theta}+\textsf{N}\right)\Theta+\left(\frac{d \textsf{N}'}{dt}\Theta'-\frac{d \textsf{N}}{dt}\Theta\right)_{t=t_0}
    + \left\langle \! \left\langle \mathbf{W},\mathbf{A}\right\rangle \!  \right\rangle  \right]dV \nonumber\\
  && ~~~~ =\int_{\mathcal{U}}\rho_0\left[\left\langle \! \left\langle\mathbf{B},\mathbf{W}\right\rangle \!\right\rangle+(R'-R)|_{t=t_0}
    \right]dV \nonumber \\
  && ~~~~~~~  +\int_{\mathcal{U}}\left[\left\langle \! \left\langle\operatorname{Div}\mathbf{P},\mathbf{W}\right\rangle \!
    \right\rangle
    +\boldsymbol{\tau}:(\nabla\mathbf{W})^{\flat}-\left(\operatorname{Div}\mathbf{Q}'-\operatorname{Div}\mathbf{Q}\right)_{t=t_0}
    \right]dV,
\end{eqnarray}
where
$\boldsymbol{\tau}:(\nabla\mathbf{W})^{\flat}=\tau^{ab}W_{a|b}$.
Assuming that $R$ and $\mathbf{Q}$ are transformed such that
\citep{MaHu1983}
\begin{equation}\label{thermal-covariance-assumptions}
    \frac{d \textsf{N}'}{dt}\Theta'-R'=\zeta\left(\frac{d \textsf{N}}{dt}\Theta-R\right)
    ~~~~~\textrm{and}~~~~~
    \mathbf{Q}'=\zeta \xi_{*}\mathbf{Q},
\end{equation}
we obtain
\begin{eqnarray}
  && \int_{\mathcal{U}}\left[ \frac{\partial \rho_0}{\partial t} +\frac{1}{2}\rho_0\operatorname{tr}\left(\frac{\partial \mathbf{G}}{\partial t}\right)-S_m\right]
  \left(\left\langle \! \left\langle\mathbf{V},\mathbf{W}\right\rangle \! \right\rangle
  +\frac{1}{2}\left\langle \! \left\langle\mathbf{W},\mathbf{W}\right\rangle \! \right\rangle\right)
    dV-\int_{\mathcal{U}}\left\langle \! \left\langle\operatorname{Div}\mathbf{P}+\rho_0\mathbf{B}-\rho_0\mathbf{A},\mathbf{W}\right\rangle \!
    \right\rangle dV \nonumber\\
  && ~+  \int_{\mathcal{U}}\left[ \left(\rho_0\frac{\partial \Psi}{\partial\mathbf{F}}:\nabla_0\mathbf{W}-\boldsymbol{\tau}:(\nabla\mathbf{W})^{\flat}\right)
    +z\rho_0\left(\frac{\partial \Psi}{\partial \Theta}+\textsf{N}\right)\Theta  \right]dV=0.
\end{eqnarray}
Note that
\begin{equation}
    \left(\nabla_0\mathbf{W}\right)^a{}_A=g^{ab}\left[(\nabla\mathbf{W})^{\flat}\right]_{bc}F^c{}_A.
\end{equation}
Therefore, arbitrariness of $\mathcal{U}$, $\mathbf{W}$, and $z$
implies that\footnote{In a previous footnote it was shown that
\begin{equation*}
    \frac{d}{dt}\Psi=\frac{\partial \Psi}{\partial \Theta}\dot{\Theta}
    +\frac{\partial \Psi}{\partial \mathbf{G}}:\frac{\partial\mathbf{G}}{\partial t}
    +\frac{\partial \Psi}{\partial\mathbf{F}}:\nabla_0\mathbf{V}
    =\frac{\partial \Psi}{\partial \Theta}\dot{\Theta}
    +\frac{\partial \Psi}{\partial \mathbf{G}}:\frac{\partial\mathbf{G}}{\partial t}
    +\frac{\partial\Psi}{\partial\mathbf{g}}:\mathfrak{L}_{\mathbf{V}}\mathbf{g}.
\end{equation*}
Thus
\begin{equation*}
    \frac{\partial \Psi}{\partial\mathbf{F}}:\nabla_0\mathbf{V}
    =\frac{\partial\Psi}{\partial\mathbf{g}}:\mathfrak{L}_{\mathbf{V}}\mathbf{g}.
\end{equation*}
This holds for an arbitray change of frame $\xi_t:\mathcal{S}\rightarrow\mathcal{S}$ as well, i.e. at $t=t_0$:
\begin{equation*}
    \frac{\partial \Psi}{\partial\mathbf{F}}:\nabla_0(\mathbf{V}+\mathbf{W})
    =\frac{\partial\Psi}{\partial\mathbf{g}}:\mathfrak{L}_{(\mathbf{V}+\mathbf{W})}\mathbf{g}.
\end{equation*}
Hence, for arbitrary $\mathbf{W}$
\begin{equation*}
    \frac{\partial \Psi}{\partial\mathbf{F}}:\nabla_0\mathbf{W}
    =\frac{\partial\Psi}{\partial\mathbf{g}}:\mathfrak{L}_{\mathbf{W}}\mathbf{g}.
\end{equation*}
Noting that \citep{MaHu1983} $\frac{\partial E}{\partial\mathbf{g}\circ\varphi}=\frac{\partial \Psi}{\partial\mathbf{g}\circ\varphi}$, this means that $\rho_0\frac{\partial \Psi}{\partial\mathbf{F}}=\mathbf{P}$ is equivalent to $2\rho_0\frac{\partial E}{\partial\mathbf{g}\circ\varphi}= \boldsymbol{\tau}$, i.e. the Doyle-Ericksen formula.
}
\begin{eqnarray}
  &&  \frac{\partial \rho_0}{\partial t} +\frac{1}{2}\rho_0\operatorname{tr}\left(\frac{\partial \mathbf{G}}{\partial t}\right)=S_m,\\
  && \operatorname{Div}\mathbf{P}+\rho_0\mathbf{B}=\rho_0\mathbf{A}, \\
  && \boldsymbol{\tau}^{\textsf{T}}=\boldsymbol{\tau}, \\
  && \rho_0\frac{\partial \Psi}{\partial\mathbf{F}}=\mathbf{P}, \\
  && \frac{\partial \Psi}{\partial \Theta}=-\textsf{N}.
\end{eqnarray}
Thus, we have proven the following proposition.

\paragraph{Proposition} Covariance of energy balance under spatial diffeomorphisms and temperature rescalings gives all the balance laws and the constitutive restrictions imposed by the Clausius-Duhem inequality.

\subsection{Covariance of the Entropy Production Inequality}

In this subsection we study the consequences of covariance of the
Clausius-Duhem inequality. Again, let us consider the
diffeomorphisms $\xi_t:\mathcal{S}\rightarrow\mathcal{S}$ and
monotonically increasing temperature rescalings
$\zeta_t:\mathbb{R}^+\rightarrow \mathbb{R}^+$. Let us postulate
that the entropy production inequality is invariant under the
simultaneous action of these two transformations, i.e.
\begin{equation}\label{entropy-new-frame}
    \frac{d}{dt}\int_{\mathcal{U}}\rho'_0 \textsf{N}'dV \geq  \int_{\mathcal{U}}\frac{\rho'_0
    R'}{\Theta'}dV+\int_{\partial\mathcal{U}}\frac{H'}{\Theta'}~dA
    +\int_{\mathcal{U}} \textsf{N}' \left[\frac{\partial \rho'_0}{\partial t}+\frac{1}{2}\rho'_0\operatorname{tr}\!\left(\frac{\partial \mathbf{G}}{\partial
    t}\right)\right] dV+\int_{\mathcal{U}} \rho'_0\frac{\partial E'}{\partial \mathbf{G}}:\frac{\partial\mathbf{G}}{\partial t}
    dV.
\end{equation}
Note that (\ref{thermal-covariance-assumptions}) implies that
\begin{equation}
    \frac{d \textsf{N}'}{dt}-\frac{R'}{\Theta'}=\frac{d \textsf{N}}{dt}-R~~~~~\textrm{and}~~~~~\frac{H'}{\Theta'}=\frac{H}{\Theta}.
\end{equation}
It can easily be shown that the inequality
(\ref{entropy-new-frame}) is identical to (\ref{Entropy}), i.e.
assuming the transformations
(\ref{thermal-covariance-assumptions}), entropy production
inequality is trivially covariant.

\subsection{Principle of Maximum Entropy Production}

In this subsection we use the so-called maximum entropy production principle to obtain a kinetic equation for $\dot{\mathbf{G}}$. This principle states that a non-equilibrium system with some possible constraints evolves in such a way to maximize its entropy production \citep{Ziegler1983,Rajagopal2004a}. This principle has found applications in many different fields of science. For a recent review see \cite{MartyushevSeleznev2006}. This principle has recently been used in growth mecahncis for obtaining kinetic equations for ``growth velocity gradient" \citep{LoretSimoes2005,AmbrosiGuana2007,Fusi2006}. Here, we use it in our geometric framework.

For a growing body, entropy production in a subbody $\mathcal{U}\subset\mathcal{B}$ is defined as
\begin{eqnarray}
 && \Gamma(\mathcal{U},t)=\frac{d}{dt}\int_{\mathcal{U}}\rho_0\textsf{N}dV-\int_{\mathcal{U}}\frac{\rho_0
    R}{\Theta}dV-\int_{\partial\mathcal{U}}\frac{H}{\Theta}~dA
    -\int_{\mathcal{U}} \textsf{N} \left[\frac{\partial \rho_0}{\partial t}+\frac{1}{2}\rho_0\operatorname{tr}\!\left(\frac{\partial \mathbf{G}}{\partial
    t}\right)\right] dV-\int_{\mathcal{U}} \rho_0\frac{\partial E}{\partial \mathbf{G}}:\frac{\partial\mathbf{G}}{\partial t}
    dV.
     \nonumber\\
 && ~~~~~=\int_{\mathcal{U}}\left[\rho_0\frac{d\textsf{N}}{dt}-\rho_0\frac{R}{\Theta}+\operatorname{Div}\!\left(\frac{\mathbf{Q}}{\Theta}\right)-\rho_0\frac{\partial E}{\partial \mathbf{G}}:\frac{\partial\mathbf{G}}{\partial t}
 \right]dV=\int_{\mathcal{U}} \frac{\Lambda}{\Theta}dV,
\end{eqnarray}
where $\Lambda=\Theta\left[\rho_0\frac{d\textsf{N}}{dt}-\rho_0\frac{R}{\Theta}+\operatorname{Div}\left(\frac{\mathbf{Q}}{\Theta}\right)-\rho_0\frac{\partial E}{\partial \mathbf{G}}:\frac{\partial\mathbf{G}}{\partial t}\right]$ is the rate of entropy production. Using energy balance, we have
\begin{equation}
    \Lambda=-\rho_0\frac{d \Psi}{dt}+\mathbf{P}:\nabla_0\mathbf{V}-\frac{1}{\Theta}\mathbf{d}\Theta\cdot\mathbf{Q}-\dot{\Theta}\textsf{N}.
\end{equation}
Note that
\begin{equation}
    -\rho_0\frac{d \Psi}{dt}=\dot{\Theta}\textsf{N}-\rho_0\frac{\partial \Psi}{\partial \mathbf{G}}:\frac{\partial \mathbf{G}}{\partial t}-\mathbf{P}:\nabla_0\mathbf{V}.
\end{equation}
Thus, we can write
\begin{equation}\label{Constraint}
    \Lambda=-\rho_0\frac{\partial \Psi}{\partial \mathbf{G}}:\dot{\mathbf{G}}-\frac{1}{\Theta}\mathbf{d}\Theta\cdot\mathbf{Q}.
\end{equation}
We now maximize $\Lambda$ with respect to $\dot{\mathbf{G}}$ under the constraint (\ref{Constraint}). Let us define
\begin{equation}
    \Phi=\Lambda+\lambda\left(\Lambda+\rho_0\frac{\partial \Psi}{\partial \mathbf{G}}:\dot{\mathbf{G}}+\frac{1}{\Theta}\mathbf{d}\Theta\cdot\mathbf{Q}\right),
\end{equation}
where $\lambda$ is a Lagrange multiplier. Maximizing $\Phi$ with respect to $\dot{\mathbf{G}}$ gives
\begin{equation}\label{LagrangeMultiplier}
    \frac{\partial \Lambda}{\partial \dot{\mathbf{G}}}=-\frac{\lambda}{\lambda+1}\rho_0\frac{\partial \Psi}{\partial \mathbf{G}}.
\end{equation}
Note that part of entropy production rate is constitutively given, i.e. $\Lambda=\bar{\Lambda}(\Theta,\mathbf{G},\dot{\mathbf{G}},\mathbf{F},\mathbf{g})-\frac{1}{\Theta}\mathbf{d}\Theta\cdot\mathbf{Q}$. As the simplest example let us assume that
\begin{equation}\label{Dissipation}
    \Lambda=\beta\operatorname{tr}\dot{\mathbf{G}}^2-\frac{1}{\Theta}\mathbf{d}\Theta\cdot\mathbf{Q}
    =\beta\dot{G}_{AB}\dot{G}_{CD}G^{AC}G^{BD}-\frac{1}{\Theta}\mathbf{d}\Theta\cdot\mathbf{Q}.
\end{equation}
Thus
\begin{equation}
    \frac{\partial \Lambda}{\partial \dot{G}_{AB}}=2\beta G^{AC}G^{BD}\dot{G}_{CD}.
\end{equation}
Or
\begin{equation}\label{entropy-rate}
    \dot{\mathbf{G}}^{\sharp}=\frac{1}{2\beta}\frac{\partial \Lambda}{\partial \dot{\mathbf{G}}}
    =-\frac{\lambda}{2\beta(\lambda+1)}\rho_0\frac{\partial \Psi}{\partial \mathbf{G}}.
\end{equation}
Using (\ref{entropy-rate}) and (\ref{Dissipation}) we can write
\begin{equation}\label{1}
    \Lambda=\frac{\lambda^2}{4\beta(\lambda+1)^2}\rho_0^2\frac{\partial \Psi}{\partial \mathbf{G}}:\frac{\partial \Psi}{\partial \mathbf{G}}
    -\frac{1}{\Theta}\mathbf{d}\Theta\cdot\mathbf{Q}.
\end{equation}
At the same time using (\ref{entropy-rate}) and (\ref{Constraint}) we have
\begin{equation}\label{2}
    \Lambda=\frac{\lambda}{2\beta(\lambda+1)}\rho_0^2\frac{\partial \Psi}{\partial \mathbf{G}}:\frac{\partial \Psi}{\partial \mathbf{G}}
    -\frac{1}{\Theta}\mathbf{d}\Theta\cdot\mathbf{Q}.
\end{equation}
Looking at (\ref{1}) and (\ref{2}) we see that $\lambda=-2$ and hence
\begin{equation}\label{G-Evolution-Entropy}
    \dot{\mathbf{G}}^{\sharp}=-\frac{1}{\beta}\rho_0\frac{\partial \Psi}{\partial \mathbf{G}}.
\end{equation}

\subsection{Isotropic Growth}

For isotropic growth, material metric has the following time
dependent form:
\begin{equation}
    \mathbf{G}(\mathbf{X},t)=e^{2\Omega(\mathbf{X},t)}\mathbf{G}_0(\mathbf{X}),
\end{equation}
i.e. a family of conformal material metrics model the growth. Thus
\begin{equation}
    \frac{\partial \mathbf{G}(\mathbf{X},t)}{\partial t}=2\frac{\partial\Omega}{\partial t}\mathbf{G}(\mathbf{X},t).
\end{equation}
Therefore, balance of mass is simplified to read
\begin{equation}
    \frac{\partial \rho_0(\mathbf{X},t)}{\partial t}+\rho_0(\mathbf{X},t)\frac{\partial\Omega(\mathbf{X},t)}{\partial t}e^{2\Omega(\mathbf{X},t)} \operatorname{tr}\mathbf{G}_0(\mathbf{X})=S_m(\mathbf{X},t).
\end{equation}
Given $\mathbf{G}=\mathbf{G}(\mathbf{X},t)$, one has the following
relation between volume elements at $t_0$ and $t$:
\begin{equation}
    dV(\mathbf{X},t)=\sqrt{\frac{\det \mathbf{G}(\mathbf{X},t)}{\det \mathbf{G}(\mathbf{X},t_0)}}~dV(\mathbf{X},t_0).
\end{equation}
Or
\begin{equation}
    dV(\mathbf{X},t)=e^{N\Omega(\mathbf{X},t)}dV_0(\mathbf{X}),
\end{equation}
where $N=\dim\mathcal{B}$. Mass form has the following representation
\begin{equation}
    \textsf{m}(\mathbf{X},t)=\rho_0(\mathbf{X},t)dV(\mathbf{X},t)=e^{N\Omega(\mathbf{X},t)}\textsf{m}_0(\mathbf{X}).
\end{equation}
Note that
$\dot{\textsf{m}}(\mathbf{X},t)=N\frac{\partial\Omega}{\partial
t}\textsf{m}(\mathbf{X},t)$. Mass of a subbody
$\mathcal{U}\subset\mathcal{B}$ will have the following
time-dependent form
\begin{equation}
    \textsf{M}_t(\mathcal{U})=\int_{\mathcal{U}}\textsf{m}(\mathbf{X},t).
\end{equation}
Hence
\begin{equation}
    \frac{d}{dt}\textsf{M}_t(\mathcal{U})=\int_{\mathcal{U}}N\frac{\partial\Omega}{\partial t}\textsf{m}(\mathbf{X},t).
\end{equation}
In isotropic growth there is no change in shape due to addition or
removal of mass. In the decomposition of deformation gradient one
has
\begin{equation}
    \mathbf{F}=\mathbf{F}^e\mathbf{F}^g~~~~~\textrm{and}~~~~~\mathbf{F}^g=\mathbbm{g}(\mathbf{X},t)\mathbf{I},
\end{equation}
where $\mathbf{I}$ is the identity map and $\mathbbm{g}$ is a
scalar field. We will discuss this decomposition in more detail is
\S3. In the sequel we will obtain those isotropic growth
distributions that are stress free. But let us first look at some simple examples of material metric evolution.

\subsection{Examples of Bulk Growth}

In this subsection we look at three examples of bulk growth and
show how analytical solutions for residual stresses can be generated for both isotropic
and non-isotropic growth.

\paragraph{Example 1 (Isotropic Growth of a Neo-Hookean Annulus):}

Let us consider a two-dimensional, incompressible neo-Hookean
material in a flat two-dimensional spatial manifold. The free
energy density of a neo-Hookean material in two dimensions has the
form
\begin{equation}
    \Psi=\Psi(\mathbf{X},\mathbf{C})=\mu(\operatorname{tr}\mathbf{C}-2),
\end{equation}
where $\mathbf{C}$ is the Cauchy-Green tensor, or equivalently,
the pull-back of the spatial metric, $C_{AB} = F^a{}_AF^b{}_B
g_{ab}$, and $\mu$ is a material constant. We assume that this
form holds for a growing isotropic material. In components
\begin{equation}
    \Psi=\mu\left(F^a{}_AF^b{}_Bg_{ab}G^{AB}-2\right).
\end{equation}
The ``2'' is of no particular significance: when the material
metric is fixed, it simply shifts the free energy by a constant.
When the material metric changes its contribution to the free
energy is proportional to the time-dependent material volume,
which, for a given growth distribution, is independent of the
spatial configuration. We ignore this term, and use
$\Psi=\Psi(\mathbf{X},\mathbf{C})=\mu\operatorname{tr}\mathbf{C}$
as our definition of the free energy.

Let us assume that initially the material has a flat annular shape
$R_1\le R \le R_2$ without any stresses. We would like to
calculate the stresses that occur in the new equilibrium
configuration after a rotationally symmetric growth,
$\Omega=\Omega(R,t)$. In polar coordinates, the spatial metric and
its inverse read
\begin{equation}
    \mathbf{g}=\left(%
\begin{array}{cc}
  g_{rr} & g_{r\theta} \\
  g_{\theta r} & g_{\theta\theta} \\
\end{array}%
\right)=\left(%
\begin{array}{cc}
  1 & 0 \\
  0 & r^2 \\
\end{array}%
\right),~~~
\mathbf{g}^{-1}=\left(%
\begin{array}{cc}
  g^{rr} & g^{r\theta} \\
  g^{\theta r} & g^{\theta\theta} \\
\end{array}%
\right)=\left(%
\begin{array}{cc}
  1 & 0 \\
  0 & 1/r^2 \\
\end{array}%
\right),
\end{equation}
and thus $\det \mathbf{g}=r^2$. The only nonzero connection
coefficients are:
$\gamma_{\theta\theta}^r=-r,~\gamma_{r\theta}^{\theta}=\gamma_{\theta
r}^{\theta}=1/r$. For the rotationally symmetric time-dependent
material metric we have
\begin{equation}
    \mathbf{G}=\left(%
\begin{array}{cc}
  G_{RR} & G_{R\Theta} \\
  G_{\Theta R} & G_{\Theta\Theta} \\
\end{array}%
\right)=e^{2\Omega(R,t)}\left(%
\begin{array}{cc}
  1 & 0 \\
  0 & R^2 \\
\end{array}%
\right),~~~
\mathbf{G}^{-1}=\left(%
\begin{array}{cc}
  G^{RR} & G^{R\Theta} \\
  G^{\Theta R} & G^{\Theta\Theta} \\
\end{array}%
\right)=e^{-2\Omega(R,t)}\left(%
\begin{array}{cc}
  1 & 0 \\
  0 & 1/R^2 \\
\end{array}%
\right),
\end{equation}
and thus, $\det \mathbf{G}=R^2e^{4\Omega(R,t)}$. The following
nonzero connection coefficients are needed in the balance of
linear momentum:
\begin{equation}
    \Gamma_{RR}^R=\Omega'(R,t),~\Gamma_{\Theta\Theta}^R=-R-R^2 \Omega'(R,t),~\Gamma_{R\Theta}^{\Theta}=\Gamma_{\Theta
    R}^{\Theta}=1/R+\Omega'(R,t),
\end{equation}
where $\Omega'(R,t)=\frac{\partial\Omega}{\partial R}$. Given
$\Omega=\Omega(R,t)$, we are looking for solutions of the form
\begin{equation}
  \varphi(R,\Theta)=(r,\theta)=(r(R,t),\Theta).
\end{equation}
Thus\footnote{If one does not consider the intrinsic metric and instead uses the standard metric of the Euclidean space, $\mathbf{F}$ has the following representation
\begin{equation*}
    \mathbf{F}=\left(%
\begin{array}{cc}
  r'(R,t) & 0 \\
  0 & \frac{r(R)}{R} \\
\end{array}%
\right)~~~\text{and}~~~
\mathbf{F}_g=\left(%
\begin{array}{cc}
  e^{\Omega(R,t)} & 0 \\
  0 & e^{\Omega(R,t)} \\
\end{array}%
\right).
\end{equation*}
Thus, $\det \mathbf{F}_e=\frac{rr'}{R}e^{-2\Omega}$ and hence $J_e=1$ is equivalent to (\ref{incomp}) as expected.
}
\begin{equation}\label{F}
    \mathbf{F}=\left(%
\begin{array}{cc}
  r'(R,t) & 0 \\
  0 & 1 \\
\end{array}%
\right),~~~\mathbf{F}^{-1}=\left(%
\begin{array}{cc}
  1/r'(R,t) & 0 \\
  0 & 1 \\
\end{array}%
\right).
\end{equation}
This gives the Jacobian as
\begin{equation}
    J=\frac{r\,r'}{Re^{2\Omega(R,t)}}.
\end{equation}
Incompressibility dictates that
\begin{equation}\label{incomp}
    rr'=R e^{2\Omega(R,t)}.
\end{equation}
This differential equation has the following solution
\begin{equation}
    r^2(R,t)=r_1^2(R,t)+\int_{R_1}^{R}2\xi e^{2\Omega(\xi,t)} d\xi.
\end{equation}
Note that $r_1(R)$ is not known a priori and will be obtained
after imposing the traction boundary conditions at $r_1$ and
$r_2$. In incompressible elasticity, $P^{aA}$ is replaced by
$P^{aA}-Jp(F^{-1})^{-A}{}_bg^{ab}$, where $p$ is an unknown scalar
field (pressure) that will be determined using the constraint
$J=1$ \citep{MaHu1983}, i.e.
\begin{equation}
    P^{aA}=2\mu F^a{}_BG^{AB}-p(R)(F^{-1})^{A}{}_bg^{ab}.
\end{equation}
Therefore, using (\ref{incomp}), we obtain the nonzero stress
components as
\begin{equation}
    P^{rR}= \frac{2\mu
    R}{r}-p(R)\frac{r}{R}e^{-2\Omega(R,t)}~~~~~\textrm{and}~~~~~
    P^{\theta\Theta}=\frac{2\mu}{R^2}e^{-2\Omega(R,t)}-\frac{p(R)}{r^2},
\end{equation}
where $p(R)$ is an unknown pressure.

Balance of linear momentum in components reads
\begin{equation}\label{bal-mom}
    P^{aA}{}|{_A}=\frac{\partial P^{aA}}{\partial X^A}+\Gamma^A_{AB}P^{aB}+P^{bA}\gamma^a_{bc}F^c{}_A=0.
\end{equation}
For the radial direction, $a=r$, we have
\begin{eqnarray}
  && P^{rA}{}|{_A}=\frac{\partial P^{rA}}{\partial X^A}+\Gamma^A_{AB}P^{rB}+P^{bA}\gamma^r_{bc}F^c{}_A   \nonumber \\
  && ~~~~~~~~= \frac{\partial P^{rR}}{\partial R}+\left(\Gamma^R_{RR}+\Gamma^{\Theta}_{\Theta
  R}\right)P^{rR}+P^{\theta\Theta}\gamma^{r}_{\theta\theta}F^{\theta}{}_{\Theta}
  \nonumber\\
  && ~~~~~~~~= \frac{\partial P^{rR}}{\partial R}+\left(\frac{1}{R}+2\Omega'(R,t)\right)P^{rR}-r P^{\theta\Theta}=0.
\end{eqnarray}
This gives
\begin{equation}
    p'(R)=\frac{2\mu R}{r^2}e^{2 \Omega(R,t)}\left[2\left(1+R \Omega'\right)-\frac{R^2}{r^2}e^{2\Omega(R,t)}-\frac{r^2}{R^2}e^{-2\Omega(R,t)}\right].
\end{equation}
Assuming that $p(R_i)=0$, we obtain
\begin{equation}
    p(R)=\int_{R_i}^{R}\frac{2\mu \xi}{r^2(\xi)}e^{2\Omega(\xi,t)}\left[2\left(1+\xi \Omega'(\xi)\right)-\frac{\xi^2}{r^2(\xi)}e^{2\Omega(\xi,t)}
    -\frac{r^2(\xi)}{\xi^2}e^{-2\Omega(\xi,t)}\right]d\xi.
\end{equation}
For $a=\theta$, balance of momentum (\ref{bal-mom}) gives
\begin{equation}
	P^{\theta A}{}|{_A}= \frac{\partial P^{\theta\Theta}}{\partial \Theta}
  +\Gamma^A_{A\Theta}P^{\theta\Theta}+P^{\theta R}\gamma^{\theta}_{rr}F^{r}{}_{R}
  +P^{\theta\Theta}\gamma^{\theta}_{\theta\theta}F^{\theta}{}_{\Theta}
  = \left(\Gamma^R_{R\Theta}+\Gamma^{\Theta}_{\Theta\Theta}\right)P^{\theta\Theta}=0.
\end{equation}
i.e. this equilibrium equation is trivially satisfied. Note that $\operatorname{tr}\left(\frac{\partial \mathbf{G}}{\partial t}\right)=4\Omega'$ and hence balance of mass reads
\begin{equation}
	\frac{\partial \rho_0(R,t)}{\partial t}+2\Omega'(R,t)\rho_0(R,t)=S_m(R,t).
\end{equation}
This differential equation can be easily solved for mass density. 

Note that
if one considerers a cylinder with kinematics assumptions $r=r(R,t),
\theta=\Theta, z=kZ$ for a constant $k$, the residual stresses
will be very similar to what was just calculated.

\paragraph{Example 2 (Anisotropic Growth of a Neo-Hookean Annulus):}

Let us consider an anisotropic growth represented by the following
material metric:
\begin{equation}
    \mathbf{G}=\left(%
\begin{array}{cc}
  e^{2\Omega(R,t)} & 0 \\
  0 & R^2e^{-2\Omega(R,t)} \\
\end{array}%
\right),~~~
\mathbf{G}^{-1}=\left(%
\begin{array}{cc}
  e^{-2\Omega(R,t)} & 0 \\
  0 & 1/R^2e^{2\Omega(R,t)} \\
\end{array}%
\right),
\end{equation}
and thus, $\det \mathbf{G}=R^2$. The following nonzero connection
coefficients are needed in the balance of linear momentum:
\begin{equation}
    \Gamma_{RR}^R=\Omega'(R,t),~\Gamma_{\Theta\Theta}^R=R e^{-4\Omega(R,t)}\left[\Omega'(R,t)-1\right],~\Gamma_{R\Theta}^{\Theta}=\Gamma_{\Theta
    R}^{\Theta}=1/R-\Omega'(R,t).
\end{equation}
Given $\Omega=\Omega(R,t)$, we are looking for solutions of the
form $\varphi(R,\Theta)=(r,\theta)=(r(R,t),\Theta)$. Thus, $\mathbf{F}$ and $\mathbf{F}^{-1}$ have the forms given in (\ref{F})and this gives the Jacobian as $J=\frac{r\,r'}{R}$. Incompressibility dictates that $rr'=R$.\footnote{In the classical formulation
\begin{equation*}
    \mathbf{F}=\left(%
\begin{array}{cc}
  r'(R,t) & 0 \\
  0 & \frac{r(R,t)}{R} \\
\end{array}%
\right)~~~\text{and}~~~
\mathbf{F}_g=\left(%
\begin{array}{cc}
  e^{\Omega(R,t)} & 0 \\
  0 & e^{-\Omega(R,t)} \\
\end{array}%
\right).
\end{equation*}
Thus, $J_e=\det \mathbf{F}_e=1$ would lead to the same incompressibility constraint $rr'=R$.} This simple differential equation has the following solution
\begin{equation}
    r(R,t)=\sqrt{R^2+C(t)}=\sqrt{R^2-R^2_1+r^2_1}.
\end{equation}
Note that $r_1(R,t)$ is not known a priori and will be obtained
after imposing the traction boundary conditions at $r_1$ and
$r_2$. Now, we get the nonzero stress components
as
\begin{equation}
    P^{rR}= 2\mu e^{-2\Omega(R,t)}\frac{R}{r(R,t)}-p(R,t)\frac{r(R,t)}{R}~~~~~\textrm{and}~~~~~
    P^{\theta\Theta}=2\mu\frac{e^{2\Omega(R,t)}}{R^2}-\frac{p(R,t)}{r^2(R,t)},
\end{equation}
where $p(R,t)$ is an unknown pressure.

Balance of linear momentum in components reads
\begin{equation}\label{bal-mom}
    P^{aA}{}|{_A}=\frac{\partial P^{aA}}{\partial X^A}+\Gamma^A_{AB}P^{aB}+P^{bA}\gamma^a_{bc}F^c{}_A=0.
\end{equation}
For the radial direction, $a=r$, we have
\begin{eqnarray}
  && P^{rA}{}|{_A}=\frac{\partial P^{rA}}{\partial X^A}+\Gamma^A_{AB}P^{rB}+P^{bA}\gamma^r_{bc}F^c{}_A   \nonumber \\
  && ~~~~~~~~= \frac{\partial P^{rR}}{\partial R}+\left(\Gamma^R_{RR}+\Gamma^{\Theta}_{\Theta
  R}\right)P^{rR}+P^{\theta\Theta}\gamma^{r}_{\theta\theta}F^{\theta}{}_{\Theta}
  \nonumber\\
  && ~~~~~~~~= \frac{\partial P^{rR}}{\partial R}+\frac{1}{R}P^{rR}-r P^{\theta\Theta}=0.
\end{eqnarray}
This gives
\begin{equation}
    p'(R,t)=\frac{2\mu R}{r^2}e^{-2\Omega(R,t)}\left[2-2R \Omega'(R,t)-\frac{r^2}{R^2}e^{4\Omega(R,t)}-\frac{R^2}{r^2}\right].
\end{equation}
Assuming that $p(R_1,t)=0$, we obtain
\begin{equation}
    p(R,t)=\int_{R_1}^{R} \frac{2\mu \xi}{r^2(\xi)}e^{-2\Omega(\xi,t)}\left[2-2\xi \Omega'(\xi,t)-\frac{r^2(\xi)}{\xi^2}e^{4\Omega(\xi,t)}-\frac{\xi^2}{r^2(\xi)}\right] d\xi.
\end{equation}
Note that $r^2=R^2+C$ and thus
\begin{equation}
    p(R,t)=\int_{R_1}^{R} \frac{2\mu \xi}{\xi^2+C}e^{-2\Omega(\xi,t)}\left[2-2\xi \Omega'(\xi,t)-\frac{\xi^2+C}{\xi^2}e^{4\Omega(\xi,t)}-\frac{\xi^2}{\xi^2+C}\right] d\xi.
\end{equation}
Assuming that $p(R_2,t)=0$, $C(t)$ can be calculated using the above
equation.

For $a=\theta$, balance of momentum (\ref{bal-mom}) gives
\begin{equation}
	P^{\theta A}{}|{_A}= \frac{\partial P^{\theta\Theta}}{\partial \Theta}
  +\Gamma^A_{A\Theta}P^{\theta\Theta}+P^{\theta R}\gamma^{\theta}_{rr}F^{r}{}_{R}
  +P^{\theta\Theta}\gamma^{\theta}_{\theta\theta}F^{\theta}{}_{\Theta}
  = \left(\Gamma^R_{R\Theta}+\Gamma^{\Theta}_{\Theta\Theta}\right)P^{\theta\Theta}=0.
\end{equation}
i.e. this equilibrium equation is trivially satisfied. It is seen
that a growth that results in only change in shape and no change
in volume can still result in residual stresses. Note that $\operatorname{tr}\left(\frac{\partial \mathbf{G}}{\partial t}\right)=0$ and hence balance of mass reads
\begin{equation}
	\frac{\partial \rho_0(R,t)}{\partial t}=S_m(R,t).
\end{equation}

\paragraph{Example 3 (Spherical Growth of a Neo-Hookean Hollow Sphere):}

Let us consider a hollow sphere with inner and outer radii $R_i$
and $R_o$ initially in a coordinate system $(R,\Theta,\Phi)$. Let
us denote the spatial coordinates by $(r,\theta,\phi)$. The
spatial metric has the following form
\begin{equation}
    \mathbf{g}=\left(%
\begin{array}{ccc}
  1 & 0 & 0 \\
  0 & r^2 & 0 \\
  0 & 0 & r^2\sin^2 \phi \\
\end{array}%
\right),
\end{equation}
with the nonzero connection coefficients
$\gamma^r_{\theta\theta}=-r,\gamma^r_{\phi\phi}=-r\sin^2\phi,\gamma^{\theta}_{r\theta}=\gamma^{\theta}_{\theta
}=1/r,\gamma^{\phi}_{r\phi}=\gamma^{\phi}_{\phi r}=1/r$. For
isotropic growth of the hollow sphere we consider the following
material metric
\begin{equation}
    \mathbf{G}=e^{2\Omega(R,t)}\left(%
\begin{array}{ccc}
  1 & 0 & 0 \\
  0 & R^2 & 0 \\
  0 & 0 & R^2\sin^2 \Phi \\
\end{array}%
\right).
\end{equation}
The nonzero connection coefficients are
\begin{equation}
    \Gamma^R_{RR}=\Omega',~~\Gamma^R_{\Theta\Theta}=-R-R^2\Omega'
    ,~\Gamma^R_{\Phi\Phi}=-\left(R+R^2\Omega'\right)\sin^2\Phi
    ,~\Gamma^R_{R\Theta}=\Gamma^R_{\Theta R}=\Gamma^R_{R\Phi}=\Gamma^R_{\Phi R}=\frac{1}{R}+\Omega'.
\end{equation}
Under this symmetric change of material metric (growth)
we look for solutions of the form $r=r(R,t),~\theta=\Theta,~\phi=\Phi$. Thus\footnote{In the classical formulation
\begin{equation*}
    \mathbf{F}=\left(%
\begin{array}{ccc}
  r'(R,t) & 0 & 0 \\
  0 & \frac{r(R,t)}{R} & 0 \\
  0 & 0 & \frac{r(R,t)}{R} \\
\end{array}%
\right)~~~\text{and}~~~
\mathbf{F}_g=\left(%
\begin{array}{ccc}
  e^{\Omega(R,t)} & 0 & 0 \\
  0 & e^{\Omega(R,t)} & 0 \\
  0 & 0 & e^{\Omega(R,t)} \\
\end{array}%
\right).
\end{equation*}
Hence, $J_e=\det \mathbf{F}_e=1$ would lead to the same incompressibility constraint. See \citet{ChenHoger2000} for more details.
}
\begin{equation}
    \mathbf{F}=\left(%
\begin{array}{ccc}
  r'(R) & 0 & 0 \\
  0 & 1 & 0 \\
  0 & 0 & 1 \\
\end{array}%
\right),
\end{equation}
and hence
\begin{equation}
    J=\frac{r^2}{R^2}e^{-3\Omega}r'.
\end{equation}
Incompressibility gives us
\begin{equation}
    r^3(R)=r_1^3(R)+\int_{R_i}^{R}3\xi^2e^{3\Omega(\xi,t)}d\xi.
\end{equation}
The only nonzero stresses are
\begin{equation}
    P^{rR}=2\mu\frac{R^2}{r^2}e^{\Omega}-p\frac{r^2}{R^2}e^{-3\Omega},~~~
    P^{\theta\Theta}=\frac{2\mu}{R^2}e^{-2\Omega}-\frac{p}{r^2},~~~
    P^{\phi\Phi}=\frac{2\mu}{R^2\sin^2\Phi}e^{-2\Omega}-\frac{p}{r^2\sin^2\Phi}.
\end{equation}
Again, the equilibrium equations for $P^{\theta\Theta}$ and
$P^{\phi\Phi}$ are trivially satisfied. The only nontrivial
equilibrium equation reads
\begin{eqnarray}
  && P^{rA}{}|{_A}=\frac{\partial P^{rA}}{\partial X^A}+\Gamma^A_{AB}P^{rB}+P^{bA}\gamma^r_{bc}F^c{}_A   \nonumber \\
  && ~~~~~~~~= \frac{\partial P^{rR}}{\partial R}+\left(\Gamma^R_{RR}+\Gamma^{\Theta}_{\Theta R}+
  \Gamma^{\Phi}_{\Phi R}\right)P^{rR}+P^{\theta\Theta}\gamma^{r}_{\theta\theta}F^{\theta}{}_{\Theta}+P^{\phi\Phi}\gamma^{r}_{\phi\phi}F^{\phi}{}_{\Phi}
  \nonumber\\
  && ~~~~~~~~= \frac{\partial P^{rR}}{\partial R}+\left(\frac{2}{R}+3\Omega'\right)P^{rR}-r P^{\theta\Theta}-r\sin^2\Phi P^{\phi\Phi}=0.
\end{eqnarray}
This gives
\begin{equation}
    p'(R,t)=\frac{4\mu R^4}{r^4}e^{4\Omega(R,t)}\left[\frac{2}{R}+\Omega'(R,t)-\frac{R^2}{r^3}e^{3\Omega(R,t)}--\frac{r^3}{R^4}e^{-3\Omega(R,t)}\right].
\end{equation}
Assuming that $p(R_1)=0$, we obtain
\begin{equation}
    p(R,t)=\int_{R_1}^{R} \frac{4\mu \xi^4}{r^4(\xi)}e^{4\Omega(\xi,t)}\left[\frac{2}{\xi}+\Omega'(\xi,t)
    -\frac{\xi^2}{r^3(\xi)}e^{3\Omega(\xi,t)}-\frac{r^3(\xi)}{\xi^4}e^{-3\Omega(\xi,t)}\right] d\xi.
\end{equation}
Note that $\operatorname{tr}\left(\frac{\partial \mathbf{G}}{\partial t}\right)=6\Omega'$ and hence balance of mass reads
\begin{equation}
	\frac{\partial \rho_0(R,t)}{\partial t}+3\Omega'(R,t)\rho_0(R,t)=S_m(R,t).
\end{equation}
This differential equation can be easily solved for mass density.

\subsection{Visualizing Material Manifolds with Evolving Metrics}

In our geometric theory, we model growth in a fixed material manifold $\mathcal{B}$. We can visualize the evolution of $\mathbf{G}(t)$ by embedding $\mathcal{B}$ in some material ambient space $\mathcal{X}$ with a fixed metric $\mathbf{H}$. For us this larger space would be the Euclidean space with its standard metric. Consider a one-parameter family of isometric embeddings $\iota_t:\mathcal{B}\hookrightarrow\mathcal{X}$, i.e. $\iota_t^*\mathbf{H}=\mathbf{G}(t)$. For the sake of simplicity, let us restrict ourselves to rotationally symmetric metrics, i.e. we look at metrics of the form
\begin{equation}
    \mathbf{G}=\left(%
\begin{array}{ccc}
  M^2(R,t) & 0  \\
  0 & N^2(R,t) \\
\end{array}%
\right),
\end{equation}
in some coordinate patch $(R,\Theta)$, i.e. metric has the form: $M^2(R,t)dR^2+N^2(R,t)d\Theta^2$, where $t$ is time. Note that $M$ and $N$ are independent of $\Theta$. We now look for solutions in the set of surfaces of revolution. Let us consider a time-dependent curve $\gamma(s,t)=(\rho(s,t),\xi(s,t))$ in the plane. The surface obtained from this curve by revolution about $z$-axis has the following parametric representation:
\begin{equation}
    \Phi(s,\Theta,t)=(\rho(s,t)\cos \Theta,\rho(s,t) \sin \Theta,\xi(s,t)).
\end{equation}
The induced Euclidean metric is \citep{Peterson1997}:
\begin{equation}
    \Phi^*\left(dX^2+dY^2+dZ^2\right)=\left(\dot{\rho}^2(s,t)+\dot{\xi}^2(s,t)\right)ds^2+\rho(s,t)^2 d\Theta^2,
\end{equation}
where a superimposed dot means differentiation with respect to $s$. Given $M(R,t)^2dR^2+N(R,t)^2d\Theta^2$, let us assume that $\rho(s,t)=N(s,t)$ and hence
\begin{equation}
    \dot{\xi}(s,t)=\sqrt{M^2(s,t)-\dot{N}^2(s,t)}.
\end{equation}
Therefore
\begin{equation}
    \xi(s,t)=\int_{s_0}^{s}\sqrt{M^2(\ell,t)-\dot{N}^2(\ell,t)}d\ell.
\end{equation}
Of course, a solution may not exist. This happens when $M^2<\dot{N}^2$. This is not surprising as not every rotationally symmetric metric arises from a surface of revolution. In the following we consider an initially stress-free annulus under different rotationally symmetric growth distributions.

\paragraph{Example 1:} Consider isotropic growth, i.e.
\begin{equation}
    \mathbf{G}=  \left(%
    \begin{array}{cc}
    e^{2\Omega(R)} & 0 \\ 
    0 & R^2e^{2\Omega(R)} \\ 
  \end{array}\right),~~~M=e^{\Omega},~~~N=Re^{\Omega}.
\end{equation}
Hence $M^2-\dot{N}^2=-R\Omega'e^{2\Omega}(R\Omega'+2)$. Let us look at two cases:
\begin{enumerate}
\item [i)] $\Omega(R)=-R$: We have $M^2-\dot{N}^2=e^{-2R}(2R-R^2)$, which for $0<R<2$ gives the material manifold shown in Fig. \ref{Isotropic}(left).
\item [ii)] $\Omega(R)=-R^2$: We have $M^2-\dot{N}^2=4R^2e^{-2R^2}(1-R^2)$, which for $0<R<1$ gives the material manifold shown in Fig. \ref{Isotropic}(right).
\begin{figure}[hbt]
\begin{center}
\includegraphics[scale=0.6,angle=0]{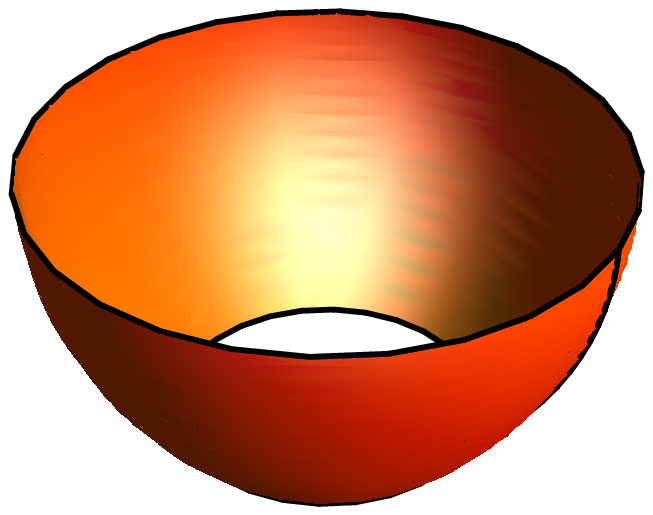}
\includegraphics[scale=0.6,angle=0]{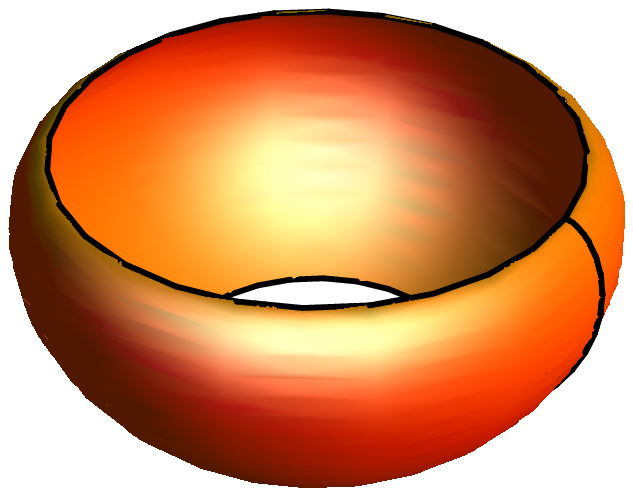}
\end{center}{\vskip -0.5 in}
\caption{\footnotesize Visualization of the material manifolds of two isotropic growth distribution of an annulus as embeddings in $\mathbb{R}^3$. Left: $\Omega(R)=-R$. Right: $\Omega(R)=-R^2$.} \label{Isotropic}
\end{figure}
\end{enumerate}

\paragraph{Example 2:} We look at anisotropic metric evolutions represented by
\begin{equation}
    \mathbf{G}=  \left(%
    \begin{array}{cc}
    e^{2\Omega(R)} & 0 \\ 
    0 & R^2e^{2\Pi(R)} \\ 
  \end{array}\right),~~~M=e^{\Omega},~~~N=Re^{\Pi}.
\end{equation}
We look at two cases:
\begin{enumerate}
\item [i)] $\Omega(R)=\cos^2 R$ and $\Pi(R)=0$: We have $M^2-\dot{N}^2=e^{2\cos^2 R}-1>0$. The material manifold shown in Fig. \ref{Anisotropic}(left).
\item [ii)] $\Omega(R)=0$ and $\Pi(R)=-\ln R^2$: We have $M^2-\dot{N}^2=1-\frac{1}{R^4}$, which for $R>1$ gives the material manifold shown in Fig. \ref{Anisotropic}(right).
\vskip 0.3 in
\begin{figure}[hbt]
\begin{center}
\includegraphics[scale=0.5,angle=0]{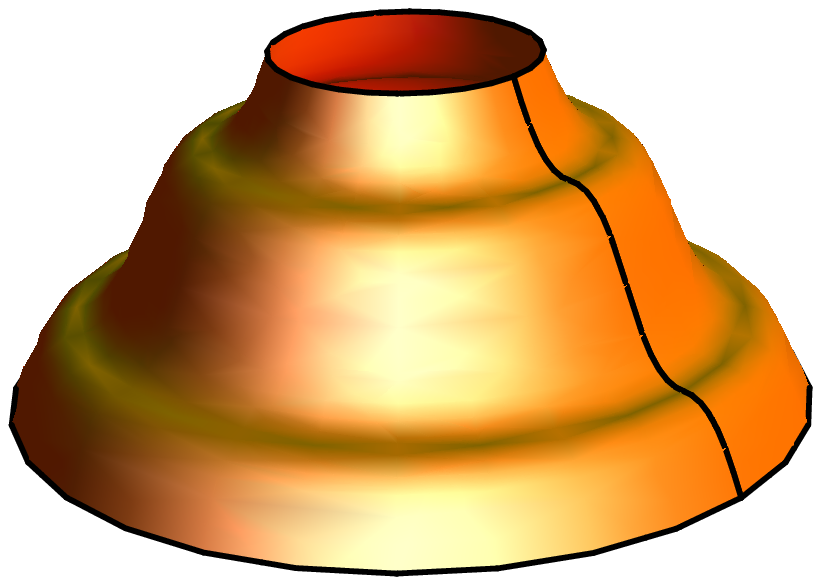}
\includegraphics[scale=0.5,angle=0]{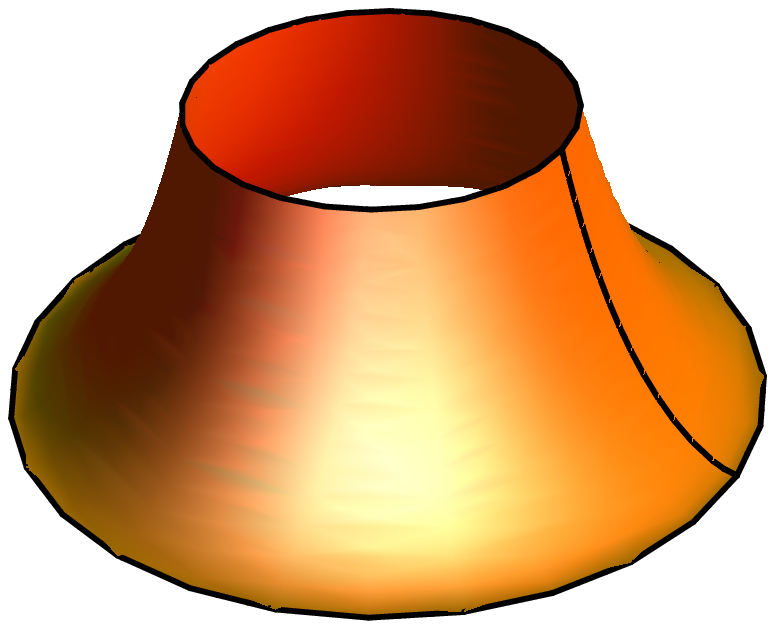}
\end{center}{\vskip -0.3 in}
\caption{\footnotesize Visualization of the material manifolds of two anisotropic growth distribution of an annulus as embeddings in $\mathbb{R}^3$. Left: $\Omega(R)=\cos^2R,~\Pi(R)=0$. Right: $\Omega(R)=0,~\Pi(R)=-\ln R^2$.} \label{Anisotropic}
\end{figure}
\end{enumerate}

\subsection{Stress-Free Isotropic Growth}

In the context of growth mechanics, \citet{TakamizawaMatsuda1990} realized that having a stress-free
configuration is equivalent to vanishing of Riemann's curvature
tensor, although they did not present any detailed calculations.
In this subsection we study this problem in detail and obtain
stress-free isotropic growth distributions in both two and three
dimensions.

Let us first review some basic concepts in Riemannian geometry.
For $\pi:E\rightarrow \mathcal{S}$ a vector bundle over a manifold
$\mathcal{S}$, $\mathcal{E}(\mathcal{S})$ the space of smooth
sections of $E$, and $\mathcal{X}(\mathcal{S})$ the space of
vector fields on $\mathcal{S}$, a connection on $E$ is a map
$\nabla:\mathcal{X}(\mathcal{S})\times\mathcal{E}(\mathcal{S})\rightarrow\mathcal{E}(\mathcal{S})$
such that $\forall~f,f_1,f_2\in
C^{\infty}(\mathcal{S}),~\forall~a_1,a_2\in\mathbb{R}$
\begin{eqnarray}
  & a)& \nabla_{f_1\mathbf{X}_1+f_2\mathbf{X}_2}\mathbf{Y}=f_1\nabla_{\mathbf{X}_1}\mathbf{Y}
        +f_2\nabla_{\mathbf{X}_2}\mathbf{Y}, \\
  & b)& \nabla_{\mathbf{X}}(a_1\mathbf{Y}_1+a_2\mathbf{Y}_2)=a_1\nabla_{\mathbf{X}}(\mathbf{Y}_1)
        +a_2\nabla_{\mathbf{X}}(\mathbf{Y}_2), \\
  & c)&
  \nabla_{\mathbf{X}}(f\mathbf{Y})=f\nabla_{\mathbf{X}}\mathbf{Y}+(\mathbf{X}f)\mathbf{Y}.
\end{eqnarray}
A linear connection on $\mathcal{S}$ is a connection on
$T\mathcal{S}$, i.e.,
$\nabla:\mathcal{X}(\mathcal{S})\times\mathcal{X}(\mathcal{S})\rightarrow\mathcal{X}(\mathcal{S})$.
In a local chart $\{x^i\}$
\begin{equation}
    \nabla_{\partial_i}\partial_j=\gamma_{ij}^k
    \partial_k,
\end{equation}
where $\gamma_{ij}^k$ are Christoffel symbols of the connection
and $\partial_i=\frac{\partial}{\partial x^i}$. A linear
connection is said to be compatible with the metric of the
manifold if
\begin{equation}
    \nabla_{\mathbf{X}}\left\langle\!\left\langle \mathbf{Y},\mathbf{Z}\right\rangle\!\right\rangle
    =\left\langle\!\left\langle \nabla_{\mathbf{X}}\mathbf{Y},\mathbf{Z}\right\rangle\!\right\rangle
    +\left\langle\!\left\langle \mathbf{Y},\nabla_{\mathbf{X}}\mathbf{Z}\right\rangle\!\right\rangle.
\end{equation}
One can show that $\nabla$ is compatible with
$\mathbf{g}$ if and only if
$\nabla\mathbf{g}=\mathbf{0}$. Torsion of a
connection is defined as
\begin{equation}
    \boldsymbol{\mathcal{T}}(\mathbf{X},\mathbf{Y})=\nabla_{\mathbf{X}}\mathbf{Y}-\nabla_{\mathbf{Y}}\mathbf{X}
    -[\mathbf{X},\mathbf{Y}],
\end{equation}
where
\begin{equation}
    [\mathbf{X},\mathbf{Y}](F)=\mathbf{X}(\mathbf{Y}(F))-\mathbf{Y}(\mathbf{X}(F))~~~~~\forall~F\in C^{\infty}(\mathcal{S}),
\end{equation}
is the commutator of $\textbf{X}$ and $\textbf{Y}$.
$\nabla$ is symmetric if it is torsion-free, i.e. $
    \nabla_{\mathbf{X}}\mathbf{Y}-\nabla_{\mathbf{Y}}\mathbf{X}
    =[\mathbf{X},\mathbf{Y}]$. According to the Fundamental Lemma of Riemannian Geometry
\citep{Lee1997} on any Riemannian manifold
$(\mathcal{S},\mathbf{g})$ there is a unique linear connection
$\nabla$, the Levi-Civita connection, that is
compatible with $\mathbf{g}$ and is torsion-free with the
following Christoffel symbols
\begin{equation}
    \gamma_{ij}^k=\frac{1}{2}g^{kl}\left(\frac{\partial g_{jl}}{\partial x^i}+\frac{\partial g_{il}}{\partial x^j}
    -\frac{\partial g_{ij}}{\partial x^l}\right).
\end{equation}

Curvature tensor $\boldsymbol{\mathcal{R}}$ of a Riemannian
manifold $(\mathcal{S},\mathbf{g})$ is a $\begin{pmatrix}  1 \\
  3 \end{pmatrix}$-tensor $\boldsymbol{\mathcal{R}}:T^*_{\mathbf{x}}\mathcal{S}\times
T_{\mathbf{x}}\mathcal{S}\times T_{\mathbf{x}}\mathcal{S}\times
T_{\mathbf{x}}\mathcal{S}\rightarrow\mathbb{R}$ defined as
\begin{equation}
    \boldsymbol{\mathcal{R}}(\alpha,\mathbf{w}_1,\mathbf{w}_2,\mathbf{w}_3)=\alpha\left(\nabla_{\mathbf{w}_1}\nabla_{\mathbf{w}_2}\mathbf{w}_3
    -\nabla_{\mathbf{w}_2}\nabla_{\mathbf{w}_1}\mathbf{w}_3
    -\nabla_{[\mathbf{w}_1,\mathbf{w}_2]}\mathbf{w}_3\right)
\end{equation}
for $\alpha\in
T^*_{\mathbf{x}}S,~\mathbf{w}_1,\mathbf{w}_2,\mathbf{w}_3\in
T_{\mathbf{x}}S$. In a coordinate chart $\{x^a\}$
\begin{equation}
    \mathcal{R}^a{}_{bcd}=\frac{\partial \gamma^a_{bd}}{\partial x^c}-\frac{\partial \gamma^a_{bc}}{\partial x^d}
    +\gamma^a_{ce}\gamma^e_{bd}-\gamma^a_{de}\gamma^e_{bc}.
\end{equation}
Note that for an arbitrary vector field $\mathbf{w}$
\begin{equation}
    w^a{}_{|bc}-w^a{}_{|cb}=\mathcal{R}^a{}_{bcd}w^d+\mathcal{T}^d{}_{cb}w^a{}_{|d}.
\end{equation}
An $n$-dimensional Riemannian manifold is flat if it is isometric
to Euclidean space. A Riemannian manifold is flat if and only if
its curvature tensor vanishes \citep{Lee1997,Spivak1999,Berger2003}. Ricci
curvature is defined as
\begin{equation}
    R_{ab}=\mathcal{R}^c{}_{acb}.
\end{equation}
The trace of Ricci curvature is called scalar curvature:
\begin{equation}
    \textsf{R}=R_{ab}g^{ab}.
\end{equation}

In dimensions two and three Ricci curvature algebraically
determines the entire curvature tensor. In dimension three
\citep{Hamilton1982}:
\begin{equation}
    \mathcal{R}_{abcd}=g_{ac}R_{bd}-g_{ad}R_{bc}-g_{bc}R_{ad}+g_{bd}R_{ac}-\frac{1}{2}\textsf{R}\left(g_{ac}g_{bd}-g_{ad}g_{bc}\right).
\end{equation}
In dimension two $R_{ab}=\textsf{R}g_{ab}$, and hence scalar
curvature completely characterizes the curvature tensor and is
twice the Gauss curvature. Note that any one-dimensional metric is flat. In the following we obtain the stress-free growth distributions in dimensions two and three.

\paragraph{i) The two-dimensional case.} Consider a
two-dimensional shell restricted to live on a flat planar surface
between two rigid planes. We assume that with no external or body
forces, initially the shell is stress-free. Can one find the
growth distributions that will result in equilibrium
configurations with zero stress? Uniform growth will obviously
result in uniform expansion/contracction, and hence no stress. Are
there other isotropic growth distributions with this property?

The spatial distances between material points are measured by the
ambient space metric (the ``spatial metric''), which is Euclidean.
A given growth distribution will result in a change in the
material metric. A configuration will be stress-free if there is
no ``stretch'' in the material, i.e., if the material distance
between two points is the same as the spatial distance. This can
happen only if the two metric tensors (spatial and material) give
the same distance measurements between nearby material points, i.e. if they are isometric. As
the spatial metric is assumed to be Euclidean, this means that the
material metric, after the change due to a given growth
distribution, must be Euclidean.

Riemann defined the curvature tensor of the metric and proved that
a metric is flat, i.e., it can be brought into the Euclidean form
$\delta_{IJ}$ locally by a coordinate transformation, if and only
if its curvature tensor is zero \citep{Lee1997,Spivak1999,Berger2003}. It
turns out that in dimension two, a weaker requirement is
sufficient \citep{Berger2003}: a metric is flat if and only if its
scalar curvature (the Ricci scalar) is zero. Let us now apply this
condition to a two-dimensional metric that is obtained from a
non-uniform growth distribution on an initially stress-free,
planar shell, i.e., $G_{IJ}=e^{2\Omega}\delta_{IJ}$. The Ricci
scalar for a metric of this form is given by \citep{Wald1984}
\begin{equation}
  \textsf{R} = -2\,e^{-2\Omega} \nabla^2\Omega\,.
\end{equation}
Thus, $\textsf{R}=0$ requires $\nabla^2\Omega=0$, i.e., $\Omega$ has to be a harmonic function. Note
that here $\nabla^2$ is the spatial Laplacian. Growth is a slow
process compared to elastic deformations and therefore time can be
treated as a parameter and hence inertial effects can be ignored. Hence, time in $\Omega$ is treated as a parameter.

It is worth emphasizing the distinction between local and global
flatness, and the implications for stress-free growth
distributions. Although the surface of a right circular
cylinder in three dimensions looks curved, it is locally,
intrinsically flat. For any given point on the cylinder, one can
find a finite-sized region containing the point, and a
single-valued coordinate patch on this region, for which the
metric has the Euclidean form. Physically, this means that for any
given point, we can cut some finite-sized piece containing the
point, and can lay the piece on a flat plane, without stretching
it. The surface of a sphere in three dimensions, on the other
hand, is intrinsically curved; it is impossible to make any
finite-sized piece of the sphere, no matter how small, to lie on a
flat plane without stretching it. The curvature condition
$\textsf{R}=0$ (or $\nabla^2 \Omega=0$) is local. Making a full
cylinder to lie in a plane nicely (i.e., without tearing, folding,
or stretching it) being impossible is due to the global topology
of the cylinder; local restrictions on curvature cannot constrain
the global properties sufficiently.

\begin{figure}[hbt]
\begin{center}
\includegraphics[scale=0.7,angle=0]{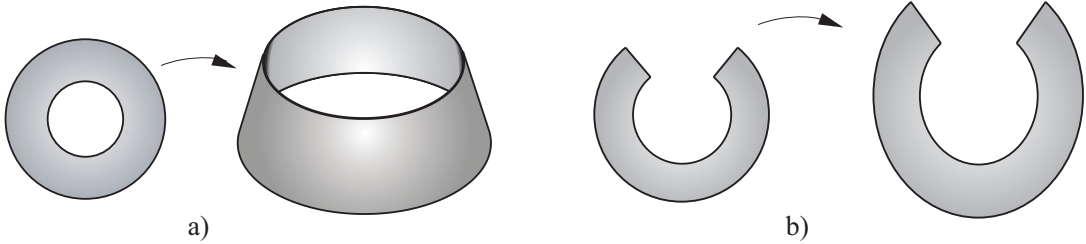}
\end{center} {\vskip -0.2 in}
\caption{\footnotesize a) Zero-stress growth of an annulus to a
cone. b) Zero-stress growth of a simply-connected piece of an
annulus.} \label{ZeroStressGrowth}
\end{figure}

Let us specialize to the case where $\Omega$ depends only on the
radial coordinate $R$ of an initially flat annular piece of a
material, $R_0\le R \le R_1$. The flatness condition gives
\begin{equation}
  \nabla^2 \Omega = \frac{1}{R}\frac{\partial}{\partial R}\left(R\,\frac{\partial \Omega(R,t)}{\partial R}\right) = 0\,.
\end{equation}
Solving this gives
\begin{equation}
  \label{radial-stress-free}
  e^{2\Omega} = \xi(t) R^{2\eta(t)},
\end{equation}
where $\xi>0$ and $\eta$ are time-dependent constants. The
metric rescaling (\ref{radial-stress-free}), with the proper
identifications, is describing an annular piece from a conical
surface, with deficit angle $\xi = 2\pi(1-1/|c|)$, where $c=\frac{1}{1+\eta}$
\citep{OzakinYavari2009}. Now, one can show that it is not
possible to make such a conical surface lie on the plane without
tearing, stretching, or folding it. Thus, starting with an annular
shell between two rigid planes, a growth distribution of the form
(\ref{radial-stress-free}) will indeed result in stresses,
although the related material metric is intrinsically flat (see
Fig. \ref{ZeroStressGrowth}.a). However, if the material consists
only of a simply-connected piece of the annulus (say, $R_1<R<R_2$,
$0<\Theta_1<\Theta<\Theta_2<2\pi$), the growth distribution
(\ref{radial-stress-free}) will just cause a stress-free expansion
of the material, between the two rigid planes. See Fig.
\ref{ZeroStressGrowth}.b.

\paragraph{A remark on conformally flat manifolds and growth mechanics.} A Riemannian manifold $(\mathcal{B},\mathbf{G})$ is conformally
flat if there exists a smooth map
$f:\mathcal{B}\rightarrow\mathbb{R}$ such that
$\mathbf{G}=f\boldsymbol{\delta}$, where $\boldsymbol{\delta}$ is
the Euclidean metric. In \emph{isothermal coordinates} the
conformally flat Riemannian metric has the following local form
\begin{equation}
    \mathbf{G}=f(\mathbf{X})\left(dX_1^2+...+dX_n^2\right).
\end{equation}
It is known that \citep{Berger2003} any two-dimensional Riemannian
manifold is conformally flat and the map $f$ is unique. A
corollary of this theorem in our theory of growth mechanics is
that given any smooth curved 2D stress-free solid, there exists a
unique growth distribution such that in the new (grown)
configuration, the 2D solid is flat and still stress free.
Equivalently, starting from a stress free flat sheet, it is always
possible to deform it to any smooth curved shape by growth without
imposing any residual stresses.

\paragraph{ii) The three-dimensional case.}

Let us next consider the three-dimensional case. In three
dimensions, a vanishing Ricci scalar is not sufficient to
guarantee local flatness. However, a three-dimensional metric is
flat if and only if its Ricci tensor vanishes \citep{Berger2003}.
The Ricci tensor $R_{IJ}$ of the metric $G_{IJ} =
e^{2\Omega}\stackrel{\tiny \circ}{G}_{IJ}$ is given in terms of
the Ricci tensor $\stackrel{\tiny \circ}{R}_{IJ}$ of
$\stackrel{\tiny \circ}{G}_{IJ}$ by the following relation
\citep{Wald1984}
\begin{equation}
  R_{IJ}=\stackrel{\tiny \circ}{R}_{IJ} - (n-2)\nabla_I\nabla_J\Omega
  - \stackrel{\tiny \circ}{G}_{IJ}\stackrel{\tiny \circ}{G}^{KL}\nabla_K\nabla_L\Omega + (n-2)\nabla_I\Omega\nabla_J\Omega-
    (n-2)\stackrel{\tiny \circ}{G}_{IJ}\stackrel{\tiny \circ}{G}^{KL}\nabla_K\Omega\nabla_L\Omega\,,
\end{equation}
where $n$ is the dimensionality. Now, once again, assume that the
initial metric $\stackrel{\tiny \circ}{G}_{IJ} = \delta_{IJ}$,
$\stackrel{\tiny \circ}{R}_{IJ}=0$, and $n=3$, and replace the
covariant derivatives with partial derivatives. This gives
\begin{equation}
  R_{IJ} = -\partial_I\partial_J\Omega -
    \delta_{IJ}\delta^{KL}\partial_K\partial_L\Omega +
    \partial_I \Omega \partial_J\Omega - \delta_{IJ}\delta^{KL}\partial_K\Omega\partial_L\Omega=0.
\end{equation}
This gives the following system of nonlinear partial differential equations
in terms of $\Omega$:
\begin{eqnarray}
  && \label{flatness-1}  \Omega_{,12}=\Omega_{,1}\Omega_{,2}, \\
  && \label{flatness-2}  \Omega_{,13}=\Omega_{,1}\Omega_{,3}, \\
  && \label{flatness-3}  \Omega_{,23}=\Omega_{,2}\Omega_{,3}, \\
  && \label{flatness-4}  \Omega_{,11}+\nabla^2\Omega+\Omega_{,2}^2+\Omega_{,3}^2=0, \\
  && \label{flatness-5}  \Omega_{,22}+\nabla^2\Omega+\Omega_{,1}^2+\Omega_{,3}^2=0, \\
  && \label{flatness-6}  \Omega_{,33}+\nabla^2\Omega+\Omega_{,1}^2+\Omega_{,2}^2=0.
\end{eqnarray}
This system of nonlinear equations were solved in
\citep{OzakinYavari2009}. The general solution is
\begin{equation}
  \Omega(X^1,X^2,X^3,t)=-\ln\left\{c_0(t)\left[(X^1)^2+(X^2)^2+(X^3)^2\right]+c_1(t)X^1+c_2(t)X^2+c_3(t)X^3+c_4(t)\right\}.
\end{equation}
In a special case if $c_1=c_2=c_3=c_4=0$, we have
\begin{equation}
  \Omega(X^1,X^2,X^3)=-\ln\left(c_0R^2\right),
\end{equation}
where $R = \sqrt{(X^1)^2 + (X^2)^2 + (X^3)^2}$. In order to
understand what this solution represents physically, let us write
the metric in polar coordinates.
\begin{equation}
	dS^2 = e^{2\Omega} \left[dR^2 + R^2 (d\Theta^2 + \sin^2\Theta d \Phi^2)\right]
       = \frac{1}{c^2R^4}\left[dR^2 + R^2 (d\Theta^2 + \sin^2\Theta d
       \Phi^2)\right].
\end{equation}
Now let us define
\begin{equation}  \label{inversion}
 \tilde{R} = \frac{1}{cR}.
\end{equation}
In terms of $\tilde{R}$, the metric becomes
\begin{equation}
  dS^2 = d\tilde{R}^2 + \tilde{R}^2(d\Theta^2 + \sin^2\Theta
  d\phi^2),
\end{equation}
which is precisely the flat Euclidean metric in three dimensions.
Thus, after the growth, the metric is still flat, but the radial
coordinate in which it is manifestly so is related to the old
radial coordinate by (\ref{inversion}) (up to a simple shift of
origin). This means that, particles at the two radii $R_1 < R_2$
move to the new radii $\tilde{R}_1 > \tilde{R}_2$, after the
growth, i.e., the material gets ``inverted''. This may not be
possible for a solid ball without tearing it apart, but it is
perfectly possible for a piece from such a ball.

If only $c_4$ is nonzero, we recover the trivial uniform growth.
If only $c_1$ is nonzero and assuming that the initial material
metric is Euclidean for the half space $X^1>0$, we have
\begin{equation}
  G_{IJ}=\frac{\lambda(t)}{(X^1)^2}\delta_{IJ},
\end{equation}
where $\lambda=1/(c_1)^2$. This shows that the material manifold is conformal to the Poincar\'e
half space.

\subsection{Lagrangian Field Theory of Growing Bodies}

In the Lagrangian formulation of nonlinear elasticity, one assumes the
existence of a Lagrangian density
\begin{equation}
    \mathcal{L}=\mathcal{L}\left(\mathbf{X},t,\mathbf{G},\varphi,\dot{\varphi},\mathbf{F},\mathbf{g}\right).
\end{equation}
Lagrangian is defined in the reference configuration as
\begin{equation}
    L=\int_{\mathcal{B}}\mathcal{L}\left(\mathbf{X},t,\mathbf{G}(\mathbf{X}),\varphi(\mathbf{X}),\dot{\varphi}(\mathbf{X}),\mathbf{F}(\mathbf{X}),\mathbf{g}(\varphi(\mathbf{X}))\right)
    dV(\mathbf{X}).
\end{equation}
In the case of a growing continuum material metric will be a dynamical
variable too. Thus, for growth of an elastic body we assume the existence of a Lagrangian density $\mathcal{L}=\mathcal{L}\left(\mathbf{X},t,\mathbf{G},\varphi,\dot{\varphi},\mathbf{F},\mathbf{g}\right)$ and write the Lagrangian as
\begin{equation}
	L=\int_{\mathcal{B}}\mathcal{L}\Big(\mathbf{X},t,\mathbf{G}(\mathbf{X},t),
     \varphi(\mathbf{X},t),\dot{\varphi}(\mathbf{X},t),\mathbf{F}(\mathbf{X},t),\mathbf{g}\circ\varphi(\mathbf{X},t)\Big)
     dV(\mathbf{X}),
\end{equation}
where $dV(\mathbf{X})=\sqrt{\det \mathbf{G}}~dX^1\wedge...\wedge
dX^n=\sqrt{\det \mathbf{G}}~d\mathbf{X}$. Having the Lagrangian,
action is defined as
\begin{equation}
    S=\int_{t_0}^{t_1}L~dt
\end{equation}
and \emph{Hamilton's Principle of least action} states that
\begin{equation}
    \delta S=\mathbf{d}\mathcal{S}\cdot\left(\delta\varphi,\delta\mathbf{G}\right)=0.
\end{equation}
The problem with this formulation is that it assumes that the 
solid is a conservative system. This is obviously not correct here as growth is a dissipative process, in general. There have been
recent works on Lagrangian formulation of dissipative systems. One
idea is to use fractional derivatives and assume that Lagrangian
is a function of some non-integer time derivatives of generalized
coordinates \citep{Riewe1997a}. It is not clear how one
can use this idea for a general field theory and even if
successful how useful that theory will be. Another way of
considering dissipation in Lagrangian mechanics is to use a
Rayleigh dissipation function \citep{Mars03}.

Assume that there exists a Rayleigh dissipation function $\mathcal{R}=\mathcal{R}(\mathbf{G},\dot{\mathbf{G}})$. For a continuum with dissipative forces $\mathbf{F}$,
the Lagrange-d'Alembert Principle states that \citep{Mars03}
\begin{equation}
    \delta\int_{t_0}^{t_1}\int_{\mathcal{B}}\mathcal{L}dV dt+\int_{t_0}^{t_1}\int_{\mathcal{B}}\mathbf{F}\cdot\delta\varphi~dt=0.
\end{equation}
Assuming the existence of a dissipation potential $\mathcal{R}$ for a growing body,
the two dissipative forces are represented as
\begin{equation}
    \mathbf{F}=-\frac{\partial\mathcal{R}}{\partial \dot{\varphi}}~~~\text{and}~~~
    \mathbf{F}_G=-\frac{\partial\mathcal{R}}{\partial \dot{\mathbf{G}}}.
\end{equation}
In this case, Lagrange-d'Alembert Principle states that
\begin{equation}
	\delta\int_{t_0}^{t_1}\int_{\mathcal{B}}\mathcal{L}\left(\mathbf{X},t,\mathbf{G}
     ,\varphi,\dot{\varphi},\mathbf{F},\mathbf{g}\circ\varphi\right)
     dVdt+\int_{t_0}^{t_1}\int_{\mathcal{B}}\left(\mathbf{F}\cdot\delta\varphi+\mathbf{F}_G\cdot\delta\mathbf{G} \right)dVdt=0.
\end{equation}
For the sake of simplicity, let us consider the two variations
separately.

\vskip 0.1 in \noindent \emph{\textbf{Case 1:}} If only
deformation mapping is varied, one has
\begin{equation}
    \delta S=\mathbf{d}S\cdot\left(\delta\varphi,\mathbf{0}\right)=0.
\end{equation}
This can be simplified to read \citep{YaMaOr2006}
\begin{equation}
    \frac{\partial\mathcal{L}}{\partial\varphi^a}-\frac{d}{d t}\frac{\partial\mathcal{L}}{\partial\dot{\varphi}^a}
    -\left(\frac{\partial\mathcal{L}}{\partial F^a{}_A}\right)_{|A}
    -\frac{\partial\mathcal{L}}{\partial F^b{}_A} F^c{}_A\gamma^b_{ac}+2\frac{\partial\mathcal{L}}{\partial g_{cd}}~g_{bd}\gamma^b_{ac}=\frac{\partial\mathcal{R}}{\partial \dot{\varphi}^a}.
\end{equation}
Or
\begin{equation}
    P_a{}^A{_{|A}}+\frac{\partial\mathcal{L}}{\partial\varphi^a}+\left(F^c{}_A P_b{}^A-J\sigma^{cd}g_{bd}\right)\gamma^b_{ac}
    =\rho_0g_{ab}A^b+\frac{\partial\mathcal{R}}{\partial \dot{\varphi}^a}.
\end{equation}

\vskip 0.1 in \noindent \emph{\textbf{Case 2:}} In material representation of classical nonlinear elasticity density is not a dynamical variable but rather it is a parameter appearing in the Lagrangian. It is through a reduction process (material to spatial) that the density ends up satisfying the continuity or advection equation (see \citet{Holm1998} for more details). Here, we should note that unlike classical nonlinear elasticity, mass density varies, in general, when material metric changes. In other words, $\delta \rho_0$ and $\delta \mathbf{G}$ are related through the nonholonomic constrant of mass balance. For a similar discussion on Lagrangian formulation of fluid mechanics in Eulerian (spatial) coordinates see \citet{Brethert1970}. We know that by definition of $S_m$
\begin{equation}
    \frac{d}{dt}\int_{\mathcal{U}}\rho_0(X,t)dV=\int_{\mathcal{U}}S_m(X,t) dV=\int_{\mathcal{U}}\stackrel{\tiny \circ}{S}_m\!(X,t) d\!\stackrel{\tiny \circ}{V},
\end{equation}
where $\stackrel{\tiny \circ}{S}_m\!(X,t)$ is mass source in the initial material manifold with volume element $d\!\stackrel{\tiny \circ}{V}$. Note that $\stackrel{\tiny \circ}{S}$ is the quantity that can be given physically. Now balance of mass can be rewritten as
\begin{equation}
    \int_{\mathcal{U}}\rho_0(X,t)dV=\int_{\mathcal{U}}\rho_0(X,t_0)dV+\int_{t_0}^{t}\int_{\mathcal{U}}\stackrel{\tiny\circ}{S}_m\!(X,t) d\!\stackrel{\tiny \circ}{V}d\tau.
\end{equation}
Now for a fixed $\stackrel{\tiny\circ}{S}_m$, let us consider mass density and material metric variation fields $\rho_0(X,t;\epsilon)$ and $\mathbf{G}(X,t;\epsilon)$. For an arbitrary $\epsilon$ the above integral mass balance reads:
\begin{equation}
    \int_{\mathcal{U}}\rho_0(X,t;\epsilon)dV_{\epsilon}=\int_{\mathcal{U}}\rho_0(X,t_0;\epsilon)dV_{\epsilon}+\int_{t_0}^{t}\int_{\mathcal{U}}\stackrel{\tiny\circ}{S}_m\!(X,\tau) d\!\stackrel{\tiny \circ}{V}d\tau.
\end{equation}
Let us take derivatives with respect to $\epsilon$ of both sides, evaluate them at $\epsilon=0$ and note that all variations vanish at $t=t_0$. This gives us
\begin{equation}
    \int_{\mathcal{U}}\left(\delta\rho_0+\frac{1}{2}\rho_0\operatorname{tr}(\delta\mathbf{G})\right)dV
    =0.
\end{equation}
As $\mathcal{U}$ is arbitrary we obtain
\begin{equation}\label{variation-constraint}
    \delta\rho_0+\frac{1}{2}\rho_0\operatorname{tr}(\delta\mathbf{G})
    =0.
\end{equation}
We can write $\mathcal{L}=\rho_0\bar{\mathcal{L}}$, where $\bar{\mathcal{L}}$ is Lagrangian density per unit mass. Hence
\begin{equation}
       \delta\int_{t_0}^{t_1}\int_{\mathcal{B}}\mathcal{L}dVdt
       =\delta\int_{t_0}^{t_1}\int_{\mathcal{B}}\rho_0\bar{\mathcal{L}}dVdt
       =\int_{t_0}^{t_1}\int_{\mathcal{B}}\left[\rho_0\delta\bar{\mathcal{L}}+\bar{\mathcal{L}}\left(\delta\rho_0+\frac{1}{2}\rho_0\operatorname{tr}(\delta\mathbf{G})\right)\right]dVdt.
\end{equation}
Thus, using (\ref{variation-constraint})
\begin{equation}
       \delta\int_{t_0}^{t_1}\int_{\mathcal{B}}\mathcal{L}dVdt=\int_{t_0}^{t_1}\int_{\mathcal{B}}\delta\mathcal{L}dVdt,
\end{equation}
where in $\delta\mathcal{L}$ mass density is assumed to be fixed. Now substituting this in the Lagrange-d'Alembert Principle and assuming that $\mathbf{G}$ varies while $\varphi$ is fixed, we obtain
\begin{equation}
	\int_{t_0}^{t_1}\int_{\mathcal{B}}\left(\frac{\partial\mathcal{L}}{\partial \mathbf{G}}-\frac{\partial\mathcal{R}}{\partial \dot{\mathbf{G}}}\right)
	:\delta \mathbf{G}dVdt=0.
\end{equation}
Therefore, the corresponding Euler-Lagrange equations are
\begin{equation}
    \frac{\partial\mathcal{L}}{\partial \mathbf{G}}=\frac{\partial\mathcal{R}}{\partial \dot{\mathbf{G}}}.
\end{equation}
Note that this is very similar to what we obtained using the principle of maximum entropy production in \S2.6.

\paragraph{Example.} Assuming that $\mathcal{R}(\mathbf{G},\dot{\mathbf{G}})$ is a quadratic function, i.e. $\mathcal{R}(\mathbf{G},\dot{\mathbf{G}})=\omega\operatorname{tr}\dot{\mathbf{G}}^2=\omega \dot{G}_{AM}\dot{G}_{BN}G^{AB}G^{MN}$, we have
\begin{equation}
	\frac{\partial\mathcal{R}}{\partial \dot{G}_{AB}}=2\omega \dot{G}_{MN}G^{AM}G^{BN}.
\end{equation}
Thus
\begin{equation}\label{G-Evolution}
    \dot{\mathbf{G}}^{\sharp}=\frac{1}{2\omega}\frac{\partial\mathcal{L}}{\partial \mathbf{G}}.
\end{equation}

Note that $\mathcal{L}=\mathcal{T}-\mathcal{V}$, where
\begin{equation}
	\mathcal{T}=\frac{1}{2}\rho_0\left\langle \! \left\langle\mathbf{V},\mathbf{V}\right\rangle \! \right\rangle
	~~~~~\text{and}~~~~~
	\mathcal{V}=\rho_0 E+\mathcal{V}_B,
\end{equation}
where $\mathcal{V}_B$ is the potential of body forces. Therefore, the evolution equation (\ref{G-Evolution}) reads
\begin{equation}
    \dot{\mathbf{G}}^{\sharp}=-\frac{1}{2\omega}\rho_0\frac{\partial E}{\partial \mathbf{G}}.
\end{equation}
Note that $E=\Psi+\textsf{N}\Theta$ and hence
\begin{equation}
    \frac{\partial E}{\partial \mathbf{G}}=\left(\frac{\partial \Psi}{\partial \mathbf{G}}+\frac{\partial \Psi}{\partial \Theta}\frac{\partial \Theta}{\partial \mathbf{G}}\right)
    +\textsf{N}\frac{\partial \Theta}{\partial \mathbf{G}}=\frac{\partial \Psi}{\partial \mathbf{G}}.
\end{equation}
Thus
\begin{equation}
    \dot{\mathbf{G}}^{\sharp}=-\frac{1}{2\omega}\rho_0\frac{\partial \Psi}{\partial \mathbf{G}},
\end{equation}
which is identical to (\ref{G-Evolution-Entropy}) if we choose $\omega=\frac{1}{2}\beta$.

\section{Connection between $\mathbf{F}=\mathbf{F}_e\mathbf{F}_g$ and the Geometric Theory}

In the literature of growth mechanics the idea of multiplicative decomposition of deformation gradient into elastic and growth parts is usually attributed to
\cite{Rodriguez1994}, although it can be seen in several earlier works like \citep{KondaurovNikitin1987, TakamizawaMatsuda1990,
Takamizawa1991}. \citet{TakamizawaMatsuda1990} and \citet{Takamizawa1991}
considered a local stress-free configuration by using a
multiplicative decomposition of deformation gradient, although
this decomposition is implicit in their presentation. They
realized that the local stress-free configurations are not unique
but a corresponding metric is unique and defined a global
stress-free configuration by equipping the original reference
configuration with this metric giving it a Riemannian structure.
Here we look at this metric and its rigorous connection with
$\mathbf{F}=\mathbf{F}_e\mathbf{F}_g$.

It should be mentioned that similar ideas were used in plasticity
and thermoelasticity before the growth mechanics applications. For
the less familiar application in thermal stresses, the idea of
decomposition of deformation gradient goes back to the works of
Stojanovi\'c and his coworkers
\citep{Stojanovic1964,Stojanovic1969}. See \citet{VujosevicLubarda2002,Lubrada2004} and
\citet{OzakinYavari2009} for a detailed discussion. These
researchers extended Kondo's
\citep{Kondo1955a,Kondo1955b,Kondo1963,Kondo1964} and Bilby's
\citep{BilbyBulloughSmith1955,BilbyGardnerStroh1956} idea of local
elastic relaxation in the continuum theory of distributed defects
to the case of thermal stresses.\footnote{Note that the idea of local elastic relaxation was first proposed in the work of \citet{Eckart1948}.}

We have posed the following question in this paper: which space, as
opposed to the Euclidean space, would be compatible with a relaxed
state of the body? We claim that the answer to this question is, a
Riemannian manifold whose metric is related to the nonuniform
growth. This metric describes the relaxed state of the material
with respect to which the strains in a given configuration should
be measured. In this framework, the constitutive relations are
given in terms of the material metric, the (Euclidean) spatial
metric, and the deformation gradient $\mathbf{F}$.

Let us consider one of the above-mentioned imaginary relaxed
pieces. Relaxation of this piece corresponds to a linear
deformation (linear, because the piece is small) denoteed by
$\mathbf{F}_g$. If this piece is deformed in some arbitrary way
after the relaxation, one can calculate the induced stresses by
using the tangent map of this deformation in the constitutive
relations. In order to calculate the stresses induced for a given
deformation of the full body, we focus our attention to one such
particular piece. The deformation gradient of the full body at
this piece $\mathbf{F}$ can be decomposed as $\mathbf{F} =
\mathbf{F}_e\mathbf{F}_g$, where, by definition, $\mathbf{F}_e =
\mathbf{F}\mathbf{F}_g^{-1}$. Thus, as far as this piece is
concerned, the deformation of the body consists of a
relaxation, followed by a linear deformation given by
$\mathbf{F}_e$. The stresses induced on this piece, for an
arbitrary deformation of the body, can be calculated by
substituting $\mathbf{F}_e$ in the constitutive relations. One
should note that $\mathbf{F}_e$ and $\mathbf{F}_g$ are not
necessarily compatible. However, as long as we have a prescription
for obtaining $\mathbf{F}_e$ and $\mathbf{F}_g$ directly for a
given deformation map $\varphi$ for the body and a growth
distribution, we can calculate the stresses by the following
procedure. Note also that if the material manifold is one
dimensional the decomposition of deformation mapping into elastic
and growth parts is always possible. This is implicitly assumed,
for example, in \citep{SenanO'ReillyTresierras2008}.

For isotropic growth, one has the following expression for
$\mathbf{F}_g$:
\begin{equation}\label{ft-definition}
    (F_g)^A{}_B=\mathbbm{g}\delta^A_B.
\end{equation}
Given this formula for $\mathbf{F}_g$, we can calculate
$\mathbf{F}_e = \mathbf{F}\mathbf{F}_g^{-1}$ for a given
deformation, and use a constitutive relation that gives the
stresses in terms of $\mathbf{F}_e$. At first glance these two
approaches seem very different, however, they are related, as we
demonstrate next. We should emphasize that the following
discussion is not restricted to isotropic growth; given any
$\mathbf{F}_g$ our arguments can be repeated.

The constitutive relations of the two approaches are formulated in
terms of different quantities: $\mathbf{G}(\mathbf{X},t)$ and
$\mathbf{F}$ on one side, and $\mathbf{F}_e =
\mathbf{F}\mathbf{F}_g^{-1}$ on the other.  Let us start with our
approach, namely, assume that a constitutive relation is given in
terms of $\mathbf{G}(\mathbf{X},t)$ and $\mathbf{F}$. This takes
the form of a scalar free energy density function that depends on
$\mathbf{G}(\mathbf{X},t)$, $\mathbf{F}$, as well as on the
spatial metric tensor $\mathbf{g}$, and possibly $\mathbf{X}$
explicitly:
\begin{equation}\label{constitutive-riemann}
  \Psi=\Psi(\mathbf{X},\Theta,\mathbf{G}(\mathbf{X},t),\mathbf{F},\mathbf{g}\circ\varphi).
\end{equation}
$\mathbf{G}$, $\mathbf{F}$, and $\mathbf{g}$ are tensors,
expressed in terms of specific bases for the material and
the ambient spaces. A change of basis changes the components of
these tensors, but $\Psi$ does not change as it is a scalar. Let
us consider a change of basis from the original coordinate basis
$\mathbf{E}_A$ of the material space, with the following property
\begin{equation}
    \langle\!\langle \mathbf{E}_{A},\mathbf{E}_{B} \rangle\!\rangle_{\mathbf{G}}=G_{AB}\,,
\end{equation}
to an orthonormal basis $\hat{\mathbf{E}}_{\hat{A}}$ that
satisfies
\begin{equation}
    \langle\!\langle \hat{\mathbf{E}}_{\hat{A}},\hat{\mathbf{E}}_{\hat{B}}
    \rangle\!\rangle_{\mathbf{G}}=\delta_{\hat{A}\hat{B}}\,.
\end{equation}
A matrix $\textsf{F}_{\hat{A}}{}^B$ represents the transformation
between the two bases:
\begin{equation}\label{orthonormal-basis}
    \hat{\mathbf{E}}_{\hat{A}}=\textsf{F}_{\hat{A}}{}^B~\mathbf{E}_{B}\,.
\end{equation}
The orthonormality condition gives
\begin{equation}\label{orthonormal-f}
  \textsf{F}_{\hat{A}}{}^{C}\textsf{F}_{\hat{B}}{}^D
  G_{CD}=\delta_{\hat{A}\hat{B}}.
\end{equation}
Any $\textsf{F}_{\hat{A}}{}^{C}$ that satisfies this equation
gives an orthonormal basis. Given such an
$\textsf{F}_{\hat{A}}{}^{C}$, we can also obtain an orthonormal
basis for the dual space by using its inverse. Defining
$\textsf{F}^{\hat{C}}{}_D$ as the transposed inverse of the matrix
$\textsf{F}_{\hat{A}}{}^B$, i.e.,  $\textsf{F}_{\hat{A}}{}^B
\textsf{F}^{\hat{A}}{}_C = \delta_C^B$ and
$\textsf{F}_{\hat{A}}{}^B \textsf{F}^{\hat{C}}{}_B =
\delta_{\hat{A}}^{\hat{C}}$, we obtain the dual orthonormal basis
$\{\hat{\mathbf{E}}^{\hat{A}}\}$ in terms of the original dual
basis $\{\mathbf{E}^A\}$ by
\begin{equation}
  \hat{\mathbf{E}}^{\hat{A}} = \textsf{F}^{\hat{A}}{}_{B}
  \mathbf{E}^B.
\end{equation}

For isotropic growth, $G_{CD} =
e^{2\Omega(\mathbf{X},t)}\delta_{CD} =
\mathbbm{g}(\mathbf{X},t)^{2}\delta_{CD}$ gives\footnote{This
means that
\begin{equation}\label{alpha-theta-relation}
    e^{2\Omega}=\mathbbm{g}^2~~~~~\textrm{or}~~~e^{\Omega}=\mathbbm{g},
\end{equation}
as $\mathbbm{g}$ is always positive. If growth is anisotropic,
having an expression for $G_{CD}$ all these arguments can be
repeated.}
\begin{equation}\label{orthonormal-f-for-t}
  \textsf{F}_{\hat{A}}{}^C = \delta_{\hat{A}}^C ~e^{-\Omega(\mathbf{X},t)} =  \delta_{\hat{A}}^C~\mathbbm{g}^{-1}(\mathbf{X},t),
\end{equation}
as a solution to (\ref{orthonormal-f}). Here, $\delta_{\hat{A}}^B$
is 1 for $A=B$, and 0, otherwise, i.e., $\delta_{\hat{1}}^1 =
\delta_{\hat{2}}^2 = \delta_{\hat{3}}^3 = 1$, etc. One should note
that (\ref{orthonormal-f}) has other solutions, as well, which we
will comment on in the sequel. Now let us write the components of the
total deformation gradient $\mathbf{F}$ in the orthonormal basis
$\{\hat{\mathbf{E}}_{\hat{A}}\}$. The components are transformed
by using $\textsf{F}$ as:
\begin{equation}
  F^{a}{}_{\hat{A}} = \textsf{F}_{\hat{A}}{}^B F^{a}{}_B\,.
\end{equation}
Using (\ref{orthonormal-f-for-t}), (\ref{alpha-theta-relation}),
and (\ref{ft-definition}), we can clearly see that the components
$F^a{}_{\hat{A}}$ are given precisely by those of $\mathbf{F}_e$,
the ``elastic part'' of the deformation gradient in
$\mathbf{F}=\mathbf{F}_e\mathbf{F}_g$ approach:
\begin{equation}
    F^{a}{}_{\hat{A}} = \textsf{F}_{\hat{A}}{}^B F^{a}{}_B= \delta_{\hat{A}}^B e^{-\Omega(\mathbf{X},t)} F^{a}{}_B
  = (\mathbbm{g}(\mathbf{X},t))^{-1} \delta_{\hat{A}}^B F^{a}{}_B= (F_g^{-1})_{A}{}^B F^{a}{}_B
  = (F_e)^a{}_{A}\,.
\end{equation}
Thus, $\mathbf{F}_e$ is the original deformation gradient, written
in terms of an orthonormal basis in the material space.\footnote{Note that given $\textsf{F}_{\hat{A}}{}^B$, material metric can be recovered as
\begin{equation*}
	G_{AB}=\textsf{F}^{\hat{C}}{}_A\textsf{F}^{\hat{D}}{}_B\delta_{\hat{C}\hat{D}}.
\end{equation*}} 
We have also shown that there is no need for a mysterious ``intermediate
configuration'' as the target space of $\mathbf{F}_g$, the latter
simply gives an orthonormal frame in the material manifold, and as
such, can be treated as a linear map from the tangent space of the
material manifold to itself.

Although a coordinate basis $\{\mathbf{E}_A = \partial / \partial
X^A\}$ is not necessarily orthonormal, one can always obtain an
orthonormal basis by applying a pointwise change of basis
$\mathsf{F}_{\hat{A}}{}^B$. Moreover, giving an orthonormal basis
in this way is equivalent to giving a metric tensor at each point;
the inner product of any two vectors can be calculated by using
their components in the orthonormal basis. We have seen above that
in the context of growth mechanics, this means that a change in
the material metric due to a growth distribution can be given in
terms of the ``growth deformation gradient'' of the local
relaxation approach. Given an orthonormal basis
$\{\hat{\mathbf{E}}_{\hat{A}}\}$, it is possible to obtain another
one, $\{\hat{\mathbf{E}}'_{\hat{A}}\}$, by using an orthogonal
transformation $\Lambda_{\hat{A}}{}^{\hat{B}}$:
\begin{equation}
   \hat{\mathbf{E}}'_{A}=\Lambda_{\hat{A}}{}^{\hat{B}}~\hat{\mathbf{E}}_{\hat{B}},
\end{equation}
where $\Lambda_{\hat{A}}{}^{\hat{B}}$ satisfies
$\Lambda_{\hat{A}}{}^{\hat{C}} \Lambda_{\hat{B}}{}^{\hat{D}}
\delta_{\hat{C}\hat{D}} = \delta_{\hat{A}\hat{B}}$. Let the
relation between the original coordinate basis
$\{\mathbf{E}_{A}\}$ and the new orthonormal basis be given by the
matrix ${\textsf{F}'}_{\hat{A}}{}^B$ as follows
\begin{equation}
    \hat{\mathbf{E}}'_{\hat{A}}=\textsf{F}'_{\hat{A}}{}^B \mathbf{E}_B.
\end{equation}
The relation between $\mathsf{F}$ and $\mathsf{F}'$ is given as
\begin{equation}\label{oldf-newf}
    \textsf{F}'_{\hat{A}}{}^B=\Lambda_{\hat{A}}{}^{\hat{C}}\textsf{F}_{\hat{C}}{}^B.
\end{equation}
Going in the opposite direction, one can see that $\mathsf{F}$ and
$\mathsf{F}'$ represent the same material metric $\mathbf{G}$, if
and only if they are related through (\ref{oldf-newf}) for some
orthogonal matrix $\Lambda_{\hat{A}}{}^{\hat{B}}$. This means that
there is an $SO(3)$ ambiguity in the choice of $\mathsf{F}$, and
hence, in that of $\mathbf{F}_g$.

Using an orthonormal basis for the material manifold, we rewrite
the constitutive relation (\ref{constitutive-riemann}) as
\begin{equation}
  \Psi=\Psi(\mathbf{X},\Theta,G_{AB} =  \delta_{AB},F^a{}_{B} = (F_e)^a{}_B ,g_{ab}).
\end{equation}
Hence, given a constitutive relation $\Psi^{\textrm{Riem}}$ in our
(Riemannian) approach, one can obtain a constitutive relation
$\Psi^{\textrm{LR}}$ in the ``local relaxation'' approach by
simply going to an orthonormal basis by (\ref{orthonormal-basis})
and (\ref{orthonormal-f}), and ignoring the constant terms $G_{AB}
= \delta_{AB}$ and $g_{ab} = \delta_{ab}$ in the functional
dependence.
\begin{equation}
  \Psi^{\textrm{LR}}(\mathbf{X},\Theta,(F_e)^a{}_B) =  \Psi^{\textrm{Riem}}\left(\mathbf{X},\Theta,G_{AB} =
  \delta_{AB},F^a{}_{B} = (F_e)^a{}_B ,g_{ab}=\delta_{ab}\right).
\end{equation}
Going in the opposite direction is also possible; starting with a
free energy function for the $\mathbf{F}=\mathbf{F}_e\mathbf{F}_g$
approach, one can derive an equivalent free energy in the
geometric approach.\footnote{A simple example can make this clearer. Let us assume that free energy density in the classical approach is $\mu\operatorname{tr}\mathbf{C}_e$. In components this reads
\begin{equation*}
	\Psi=\mu(C_e)_{\hat{A}\hat{B}}\delta^{\hat{A}\hat{B}}
	=\mu\left(\textsf{F}_{\hat{A}}{}^A F^a{}_A\right)\left(\textsf{F}_{\hat{B}}{}^B F^b{}_B\right)\delta_{ab}\delta^{\hat{A}\hat{B}}
	=\mu\left(F^a{}_AF^b{}_B\delta_{ab}\right)\left(\textsf{F}_{\hat{A}}{}^A \textsf{F}_{\hat{B}}{}^B \delta^{\hat{A}\hat{B}}\right)
        =\mu F^a{}_AF^b{}_B\delta_{ab}G^{AB}.
\end{equation*}
Thus, $\Psi=\mu\operatorname{tr}_{\mathbf{G}}\mathbf{C}$.
}

\paragraph{Balance of Mass.} In the $\mathbf{F}=\mathbf{F}_e\mathbf{F}_g$ approach mass balance reads $S_m=\frac{\partial \rho_0}{\partial t}+\rho_0\operatorname{tr}\mathbf{L}_g$, where $\mathbf{L}_g=\dot{\mathbf{F}}_g\mathbf{F}_g^-1$. Usually, it is assumed that growth is density preserving \citep{LubardaHoger2002}. We show that the term $\operatorname{tr}\mathbf{L}_g$ is equivalent to $\frac{1}{2}\operatorname{tr}_{\mathbf{G}}\left(\frac{\partial \mathbf{G}}{\partial t}\right)$, where by $\operatorname{tr}_{\mathbf{G}}$ we emphasize the $\mathbf{G}$-dependence of the trace operator. Note that
\begin{equation*}
   \operatorname{tr}\left(\frac{\partial \mathbf{G}}{\partial t}\right)=\frac{\partial G_{AB}}{\partial t}G^{AB}
   =\frac{\partial}{\partial t}\left(\textsf{F}^{\hat{A}}{}_A \textsf{F}^{\hat{B}}{}_B \delta_{\hat{A}\hat{B}}\right)\left(\textsf{F}_{\hat{C}}{}^A \textsf{F}_{\hat{D}}{}^B \delta^{\hat{C}\hat{D}}\right)
   =2\dot{\textsf{F}}^{\hat{A}}{}_A \textsf{F}_{\hat{A}}{}^A=2\operatorname{tr}\mathbf{L}_g.
\end{equation*}

\paragraph{Incompressibility.} In the $\mathbf{F}=\mathbf{F}_e\mathbf{F}_g$ approach incompressibility is equivalent to $J_e=1$. In the geometric theory incompressibility means $J=1$. These are equivalent as is shown below:
\begin{equation}
  1=J=\sqrt{\frac{\det \mathbf{g}}{\det\mathbf{G}}}\det\mathbf{F}
  =\frac{1}{\sqrt{\det \left( \textsf{F}^{\hat{A}}{}_A \textsf{F}^{\hat{B}}{}_B~ \delta_{\hat{A}\hat{B}} \right)}}
  \det\left(F^a{}_{\hat{C}}\textsf{F}^{\hat{C}}{}_C\right)=\det F^a{}_{\hat{C}}=J_e.
\end{equation}

\paragraph{Absolutely parallelizable manifolds and their connection with growth
mechanics.} Whenever deformation is coupled with other phenomena,
e.g. plasticity, growth/remodeling, thermal expansion/contraction,
etc. all one can hope for is to locally decouple the elastic
deformations from the inelastic deformations. Many related works
start from a decomposition of deformation gradient
$\mathbf{F}=\mathbf{F}_e\mathbf{F}_a$, where $\mathbf{F}_e$ is the
elastic deformation gradient and $\mathbf{F}_a$ is the remaining
local deformation or anelastic deformation gradient. Given an (inelastic) growth deformation gradient, a vector in the
tangent space of $\mathbf{X}\in\mathcal{B}$, i.e. $\mathbf{W}\in
T_{\mathbf{X}}\mathcal{B}$ is mapped to another vector
$\hat{\mathbf{W}}=\mathbf{F}_a\mathbf{W}$. Traditionally, these
vectors are assumed to lie in the tangent bundle of an
``intermediate configuration." In the literature, intermediate
configuration is not clearly defined and at first glance it seems
to be more or less mysterious. These are closely related to
parallelizable manifolds (or absolutely parallelizable (AP)
manifolds)
\citep{Eisenhart1926,Eisenhart1927,YoussefSid-Ahmed2007,Wanas2008}.
In an $n$-dimensional AP-manifold $M$, one starts with a field of
$n$ linearly independent vectors $\left\{\mathbf{E}_{(A)}
\right\}$ that span the tangent vector at each point. We denote
the components of $\mathbf{E}_{(A)}$ by $\mathbf{E}_{(A)}^I$. The
dual vectors, i.e. the corresponding basis vectors for the
cotangent space are denoted by $\left\{\mathbf{E}^{(A)} \right\}$
with components $\left\{\mathbf{E}^{(A)}_I \right\}$. Note that
\begin{equation}
    \mathbf{E}^{(A)}_I\mathbf{E}_{(B)}^I=\delta^A_B~~~~~\textrm{and}~~~~~\mathbf{E}^{(A)}_I\mathbf{E}_{(A)}^J=\delta^I_J.
\end{equation}
One can equip $M$ with a connection $\Gamma^I_{JK}$ such that the
basis vectors $\left\{\mathbf{E}_{(A)} \right\}$ are covariantly
constant, i.e.\footnote{Equivalently, the tangent bundle is a
trivial bundle, so that the associated principal bundle of linear
frames has a section on $M$.}
\begin{equation}\label{parallel}
    \mathbf{E}_{(A)}^I{}_{|J}=0.
\end{equation}
Note that
\begin{equation}
    \mathbf{E}_{(A)}^I{}_{|JK}-\mathbf{E}_{(A)}^I{}_{|KJ}=\mathcal{R}^I{}_{LJK}\mathbf{E}_{(A)}^L+\mathcal{T}^L{}_{KJ}\mathbf{E}_{(A)}^I{}_{|L}.
\end{equation}
Therefore, (\ref{parallel}) implies that
\begin{equation}
    \mathcal{R}^I{}_{LJK}=0,
\end{equation}
i.e., $M$ is flat with respect to the connection $\Gamma^I_{JK}$.
Note that
\begin{equation}
    \mathbf{E}_{(A)}^I{}_{|J}=\frac{\partial \mathbf{E}_{(A)}^I}{\partial
    X^J}+\Gamma^I_{JK}\mathbf{E}_{(A)}^K.
\end{equation}
Thus
\begin{equation}
    \mathbf{E}^{(A)}_L\frac{\partial \mathbf{E}_{(A)}^I}{\partial X^J}+\Gamma^I_{LK}=0.
\end{equation}
Hence
\begin{equation}
    \Gamma^I_{JK}=-\mathbf{E}^{(A)}_J\frac{\partial \mathbf{E}_{(A)}^I}{\partial X^K}=\mathbf{E}^{(A)}_I\frac{\partial \mathbf{E}^{(A)}_J}{\partial X^K}.
\end{equation}
This connection has been used by many authors, e.g. by
\citet{BilbyBulloughSmith1955} and \citet{Kondo1955a} for
dislocations, by \citet{EpsteinElzanowski2007} for material
inhomogeneities, and by \citet{Stojanovic1964} for thermal
stresses. This connection is curvature-free by construction but has a non-vanishing torsion.

For a growing body, in the local charts $\left\{X^A\right\}$ and $\left\{U^I\right\}$
for the reference and intermediate configurations, we have
\begin{equation}
    dU^I=\left(F_g\right)^I{}_A~ dX^A.
\end{equation}
$\left(F_g\right)^I{}_A$ can be identified with $\mathbf{E}_{(A)}^I$, and hence
\begin{equation}
    \Gamma^I_{JK}=\left(F_g\right)^I{}_A\frac{\partial \left(F_g^{-1}\right)^A{}_J}{\partial
    X^K}.
\end{equation}
Note that this (growth) connection is curvature free but has a
non-vanishing torsion. In plasticity it is shown that torsion of this connection has a physical meaning; it can be identified with the dislocation density tensor. For a growing body such a quantity does not seem to have a physical interpretation and we prefer to work with a Riemannian material manifold whose curvature quantifies the tendency of the growth distribution in causing residual stresses.

\vskip 0.1in \noindent 
In summary, our geometric approach has a concrete connection with
that of $\mathbf{F}=\mathbf{F}_e\mathbf{F}_g$: in the geometric
approach we use a Riemannian manifold with a time-dependent metric
as the material manifold while
$\mathbf{F}=\mathbf{F}_e\mathbf{F}_g$ implicitly uses the same
metric but in an absolutely parallelizable manifold that is not
Riemannian. We believe that our approach is more straightforward
as we do not introduce an unnecessary torsion in the material
manifold but of course the Riemannian material manifold has a non-vanishing curvature tensor, in general.

\section{Linearized Theory of Growth Mechanics}

Geometric linearization of elasticity was first introduced by
\citet{MaHu1983} and was further developed by
\citet{YavariOzakin2008}. See also \citet{MazzucatoRachele2006}.
In this section, we start with a body with a time-dependent
material manifold and its motion in an ambient space, which is
assumed to be Euclidean. Suppose a given body with a material
metric $\mathbf{G}$ is in a static equilibrium configuration,
$\varphi$. The balance of linear momentum for this material body
reads\footnote{Growth is a ``slow" process compared to elastic
deformations and hence inertial effects can be ignored. Throughout
this paper, time is teated as a parameter.}
\begin{equation}\label{linear-momentum}
    \operatorname{Div}\mathbf{P}+\rho_0\mathbf{B}=\mathbf{0}\,.
\end{equation}
Now suppose the body grows by a small amount represented by a
small change in the material metric $\delta \mathbf{G}$. $\varphi$
will no longer describe a static equilibrium configuration. Stress
in this new equilibrium configuration $\varphi' = \varphi + \delta \varphi$ will
be $\mathbf{P}' = \mathbf{P}+\delta\mathbf{P}$. We are interested
in calculating the change in the stress (or the configuration),
for a given small amount of growth.

The linearization procedure can be formulated rigorously if
instead of thinking about two nearby configurations and the
differences between various quantities for these configurations,
we describe the situation in terms of a one-parameter family of
configurations around a reference motion, and calculate the
derivatives of various quantities with respect to the parameter.
Let $\mathbf{G}_{\epsilon}(\mathbf{X})$ be a one-parameter family
of material metrics, $\varphi_{\epsilon}$ be the corresponding
equilibrium configurations, and $\mathbf{P}_{\epsilon}$ be the
corresponding stresses. Let $\epsilon=0$ describe the reference
configuration. Now, for a fixed point $\mathbf{X}$ in the material
manifold, $\varphi_{\epsilon}(\mathbf{X})$ describes a curve in
the spatial manifold, and its derivative at $\epsilon=0$ gives a
vector $\mathbf{U(X)}$ at $\varphi(\mathbf{X})$
\citep{YavariOzakin2008}:
\begin{equation}
  \mathbf{U(X)} = \frac{d\varphi_{\epsilon}(\mathbf{X})}{d
  \epsilon}\Big|_{\epsilon=0}.
\end{equation}
Considering $\delta \varphi \approx \epsilon
\frac{d\varphi_{\epsilon}}{d \epsilon}$, we see that a more
rigorous version of $\delta \varphi$ is the vector field
$\mathbf{U}$. $\mathbf{U}$ is the geometric analogue of what is
called displacement field in classical linear elasticity.

First variation (or linearization) of
deformation gradient is defined as
\begin{equation}
    \mathcal{L}(\mathbf{F}):=\nabla_{\frac{\partial}{\partial \epsilon}}\mathbf{F}_{\epsilon}\Big|_{\epsilon=0}
    =\nabla_{\frac{\partial}{\partial \epsilon}}\left(\frac{\partial \varphi_{t,\epsilon}}{\partial \mathbf{X}}\right)\Bigg|_{\epsilon=0}
    =\nabla\mathbf{U}.
\end{equation}
Or in components
\begin{equation}\label{F-linearization}
    \mathcal{L}(\mathbf{F})^a{}_{A}=U^a{}_{|A}=\frac{\partial U^a}{\partial
    X^A}+\gamma^a_{bc}F^b{}_AU^c,
\end{equation}
where $\gamma^a_{bc}$ are the connection coefficients of the
Riemannian manifold $(\mathcal{S},\mathbf{g})$. Note that for
different values of $\epsilon$ the spatial leg of $\mathbf{F}_{\epsilon}$ lies in
different tangent spaces and this is why covariant derivative with
respect to $\frac{\partial}{\partial \epsilon}$ should be used. The right
Cauchy-Green strain tensor for the perturbed motion
$\varphi_{t,\epsilon}$ is defined as
\begin{equation}
    C_{AB}(\epsilon)=F^a{}_A(\epsilon)F^b{}_B(\epsilon)g_{ab}(\epsilon).
\end{equation}
Note that $\mathbf{C}_{\epsilon}$ lies in the same linear space for all
$\epsilon\in I$, and the first variation of $\mathbf{C}$ can be
calculated as
\begin{equation}
    \frac{d}{d \epsilon}C_{AB}(\epsilon)=\nabla_{\frac{\partial}{\partial \epsilon}}F^a{}_A(\epsilon)F^b{}_B(\epsilon)g_{ab}(\epsilon)
    +F^a{}_A(\epsilon)\nabla_{\frac{\partial}{\partial \epsilon}}F^b{}_B(\epsilon)g_{ab}(\epsilon).
\end{equation}
Therefore
\begin{equation}
    \mathcal{L}(\mathbf{C})_{AB}:=\frac{d}{d \epsilon}\Big|_{\epsilon=0}C_{AB}(\epsilon)=F^b{}_B~g_{ab}~U^a{}_{|A}+F^a{}_A~g_{ab}~U^b{}_{|B}.
\end{equation}
Transpose of the deformation gradient has the following
linearization \citep{YavariOzakin2008}: $\mathcal{L}\left(\mathbf{F}^{\textsf{T}}\right)=(\nabla\mathbf{U})^{\textsf{T}}$. Spatial and material strain tensors are defined, respectively, as
\citep{MaHu1983}
\begin{equation}
    \mathbf{e}=\frac{1}{2}(\mathbf{g}-\varphi_{t*}\mathbf{G})~~~\text{and}~~~
    \mathbf{E}=\frac{1}{2}(\varphi_{t}^*\mathbf{g}-\mathbf{G}).
\end{equation}
Or in components
\begin{equation}
    e_{ab}=\frac{1}{2}\left(g_{ab}-G_{AB}F^{-A}{}_aF^{-B}{}_b\right),~~~
    E_{AB}=\frac{1}{2}(C_{AB}-G_{AB}).
\end{equation}
We now show that linearization of $\mathbf{E}$ is related to
$\boldsymbol{\epsilon}=\frac{1}{2}\mathfrak{L}_{\mathbf{u}}\mathbf{g}$,
where $\mathbf{u}=\mathbf{U}\circ\varphi^{-1}$. We know that
\begin{equation}
	\mathcal{L}(\mathbf{C})_{AB} = g_{ab}F^a{}_AF^c{}_B~u^b{}_{|c}
    +g_{ab}F^b{}_BF^c{}_A~u^a{}_{|c} = F^a{}_AF^c{}_B~u_a{}_{|c}
    +F^b{}_BF^c{}_A~u_b{}_{|c} = 2F^a{}_AF^b{}_B~ \epsilon_{ab},
\end{equation}
where $\epsilon_{ab}=\frac{1}{2}(u_a{}_{|b}+u_b{}_{|a})$ is the
linearized strain. Therefore
\begin{equation}
    \mathcal{L}(\mathbf{C})=2\varphi_t^*\boldsymbol{\epsilon}.
\end{equation}
Thus
\begin{equation}
    \boldsymbol{\epsilon}=\varphi_{t*}\mathcal{L}(\mathbf{E}).
\end{equation}
In other words, linearized strain is the push-forward of the
linearized Lagrangian strain. Obviously, if the ambient space is
Euclidean and the coordinates are Cartesian the covariant
derivatives reduce to partial derivatives and one recovers the
classical definition of linear strain in terms of partial
derivatives, i.e.
\begin{equation}
    \epsilon_{ab}=\frac{1}{2}\left(\frac{\partial u_a}{\partial x^b}+\frac{\partial u_b}{\partial x^a}\right).
\end{equation}
Note that when the linearized strain is zero the variation field
is a Killing vector field for the spatial metric $\mathbf{g}$. In
other words, this shows that this definition of linearized strain
is consistent when the variation field generates an isometry of
the ambient space.

For the one-parameter family of material metrics $\mathbf{G}_{\epsilon}$, variation of the material metric is defined as
\begin{equation}
  \delta\mathbf{G}\approx \epsilon \, \frac{d}{d
  \epsilon}\Big|_{\epsilon=0}\mathbf{G}_{\epsilon}.
\end{equation}
In the case of isotropic growth
\begin{equation}\label{deriv-of-g}
 \frac{d}{d
  \epsilon}\Big|_{\epsilon=0}\mathbf{G}_{\epsilon}=\frac{d}{d \epsilon}e^{2\Omega_{\epsilon}}\mathbf{G}_0
  =2\frac{d \Omega_{\epsilon}}{d \epsilon}\Big|_{\epsilon=0}\mathbf{G}=\beta
  \mathbf{G},
\end{equation}
where $\beta=2\delta \Omega$. Now consider, in the absence of body
forces, the equilibrium equations
$\operatorname{Div}\mathbf{P}=\mathbf{0}$ for the family of
material metrics parametrized by $\epsilon:~\operatorname{Div}_{\epsilon}\mathbf{P}_{\epsilon}=\mathbf{0}$. Linearization of equilibrium equations is defined as
\citep{YavariOzakin2008}:
\begin{equation}\label{perturbed-linear-momentum1}
    \frac{d}{d\epsilon}\Big|_{\epsilon=0}\left(\operatorname{Div}_{\epsilon}\mathbf{P}_{\epsilon}\right)=\mathbf{0}.
\end{equation}
Once again, one should note that since the equilibrium
configuration is different for each $\epsilon$,
$\mathbf{P}_{\epsilon}$ is based at different points in the
ambient space for different values of $\epsilon$, and in order to
calculate the derivative with respect to $\epsilon$, one in
general needs to use the connection (parallel transport) in the
ambient space. For the case of Euclidean ambient space that we are
considering and a Cartesian coordinate system $\{x^a\}$,
(\ref{perturbed-linear-momentum1}) is simplified and in components
reads
\begin{equation}
    \frac{\partial P^{aA}(\epsilon)}{\partial X^A}+\Gamma^A_{AB}(\epsilon)P^{aB}(\epsilon)=0.
\end{equation}
Thus, the linearized balance of linear momentum can be written as
\begin{equation}\label{linearized-balance-mom}
    \frac{\partial}{\partial X^A}\frac{d}{d\epsilon}\Big|_{\epsilon=0}P^{aA}(\epsilon)
    +\left[\frac{d}{d\epsilon}\Big|_{\epsilon=0}\Gamma^A_{AB}(\epsilon)\right]P^{aB}+
    \Gamma^A_{AB}\frac{d}{d\epsilon}\Big|_{\epsilon=0}P^{aB}(\epsilon)=0.
\end{equation}
Note that
\begin{equation}
  P^{aA} = g^{ac}\frac{\partial\Psi}{\partial F^c_A}\,,
\end{equation}
where
$\Psi=\Psi(\mathbf{X},\Theta,\mathbf{F},\mathbf{G},\mathbf{g})$ is
the material free energy density. In calculating
$\frac{dP^{aA}(\epsilon)}{d\epsilon}$, we need to consider the
changes in both $\mathbf{F}$ and $\mathbf{G}$:
\begin{equation}
  \frac{dP^{aA}(\epsilon)}{d\epsilon} = \frac{\partial
  P^{aA}}{\partial F^{b}{}_B}\frac{d F^b{}_B}{d\epsilon} +  \frac{\partial
    P^{aA}}{\partial G_{CD}}\frac{d G_{CD}}{d\epsilon}.
\end{equation}
Let us define
\begin{equation}
  \mathbbm{A}^{aA}{}_b{}^B = \frac{\partial P^{aA}}{\partial F^b{}_B}
  =
  g^{ac}\frac{\partial^2 \Psi}{\partial F^b{}_B \partial F^c{}_A}
  ~~~~~~~\textrm{and}~~~~~~~
  \mathbbm{B}^{aACD} = \frac{P^{aA}}{G_{CD}} =
  g^{ac}\frac{\partial^2\Psi}{\partial G_{CD} \partial F^c{}_A},
\end{equation}
where the derivatives are to be evaluated at the reference
motion $\epsilon = 0$. Noting that for the case of an Euclidean ambient space (see (\ref{F-linearization}))
\begin{equation}
  \frac{d F^a{}_A}{d\epsilon}\Big|_{\epsilon=0} = \frac{\partial U^a}{\partial X^{A}}
\end{equation}
we obtain
\begin{equation}
    \frac{d}{d\epsilon}\Big|_{\epsilon=0}P^{aA}(\epsilon)
    =\mathbbm{A}^{aA}{}_b{}^BU^b{}_{,B}+\mathbbm{B}^{aACD}\delta G_{CD}.
\end{equation}
Using
\begin{equation}
    \Gamma^A_{BC}=\frac{1}{2}G^{AD}\left(\frac{\partial G_{BD}}{\partial X^C}+
    \frac{\partial G_{CD}}{\partial X^B}-\frac{\partial G_{BC}}{\partial X^D}\right)
\end{equation}
and
\begin{equation}
  \frac{d G^{AB}}{d \epsilon} =
  -G^{AC}G^{BD}\frac{d G_{CD}}{d\epsilon},
\end{equation}
we obtain
\begin{equation}
  \delta\Gamma^A_{AB}=\frac{d}{d\epsilon}\Big|_{\epsilon=0}\Gamma^A_{AB}(\epsilon)=-G^{CD}\delta
  G_{CD}\Gamma^A_{AB}+\frac{1}{2}G^{AD}\left[\frac{\partial \delta G_{BD}}{\partial X^C}
  +\frac{\partial \delta G_{CD}}{\partial X^B}-\frac{\partial \delta G_{BC}}{\partial
  X^D}\right].
\end{equation}
In the case of isotropic growth, this is reduced to
\begin{equation}
    \frac{d}{d\epsilon}\Big|_{\epsilon=0}\Gamma^A_{AB}(\epsilon)=\frac{3}{2}\frac{\partial \beta}{\partial X^B}.
\end{equation}
With these results, the linearized balance of linear momentum
(\ref{perturbed-linear-momentum1}) becomes
\begin{equation}
    \left(\mathbbm{A}^{aA}{}_b{}^BU^b{}_{,B}\right)_{,A}+\left(\mathbbm{B}^{aACD}\delta G_{CD}\right)_{,A}+\frac{3}{2}\frac{\partial \beta}{\partial X^B}P^{aB}=0.
\end{equation}
Assuming that $\boldsymbol{\mathbbm{A}}$ and
$\boldsymbol{\mathbbm{B}}$ are independent of $\mathbf{X}$, the
linearized equilibrium equations are simplified to read
\begin{equation}
    \mathbbm{A}^{aA}{}_b{}^B \frac{\partial^2 U^b}{\partial X^A\partial X^B}
    +\mathbbm{B}^{aACD}G_{CD} \frac{\partial \beta}{\partial X^A}+\frac{3}{2}\frac{\partial \beta}{\partial X^B}P^{aB}=0.
\end{equation}
If the initial configuration is stress-free, we have
\begin{equation}\label{linearized-momentum}
    \mathbbm{A}^{aA}{}_b{}^B \frac{\partial^2 U^b}{\partial X^A\partial X^B}
    =-\mathbbm{B}^{aACD}G_{CD} \frac{\partial \beta}{\partial X^A}.
\end{equation}
Let us now simplify the above linearized equations for a specific class of elastic materials.

\paragraph{Saint-Venant-Kirchhoff materials.}
Saint-Venant-Kirchhoff materials have a constitutive relation that
is analogous to the linear isortropic materials, namely, the
second Piola-Kirchhoff stress $\mathbf{S}$ is given in terms of
the Lagrangian strain
$\mathbf{E}=\frac{1}{2}(\mathbf{C}-\mathbf{G})$ as
\citep{MaHu1983} $\mathbf{S}=\lambda (\operatorname{tr}\mathbf{E})\mathbf{G}^{-1}+2\mu \mathbf{E}$ or in components
\begin{equation}\label{Saint-Venant-Kirchhoff}
    S^{CD}=\lambda E_{AB}G^{AB}G^{CD}+2\mu E^{CD}
    =\frac{\lambda}{2}(C_{AB}G^{AB}-3)G^{CD}+\mu(C_{AB}G^{AC}G^{BD}-G^{CD}),
\end{equation}
where $\lambda=\lambda(\mathbf{X})$ and $\mu=\mu(\mathbf{X})$ are
two scalars characterizing the material properties. We can obtain
the tensor $\mathbbm{B}^{aCAB}$ from $\mathbf{S}$ as follows
\begin{equation}
    \mathbbm{B}^{aCAB}=\frac{\partial}{\partial G_{AB}}\left(g^{ab}\frac{\partial \psi}{\partial F^b{}_C}\right)
    =\frac{\partial P^{aC}}{\partial G_{AB}}=F^a{}_D\frac{\partial S^{CD}}{\partial G_{AB}}.
\end{equation}
Using
\begin{equation}
    \frac{\partial G^{AB}}{\partial G_{MN}}=-G^{AM}G^{BN}
\end{equation}
we obtain
\begin{equation}
    \mathbbm{B}^{aACD}G_{CD}=-2C_{MN}F^a{}_B\left(\lambda G^{AB}G^{MN}+2\mu G^{AM}G^{BN}\right)+(3\lambda+2\mu)F^a{}_BG^{AB}.
\end{equation}
The initial metric is Euclidean; in Cartesian coordinates,
$G_{AB}=\delta_{AB}$. Since the ambient space is also Euclidean,
we can choose a Cartesian coordinate system whose axes coincide
with the initial location of the material points along the
material Cartesian axis. This will give, $F^a{}_A=\delta^a_A$,
where $a$ and $A$ both range over $1, 2, 3$. Hence
\begin{equation}
    \mathbbm{B}^{aACD}G_{CD}=-\frac{3\lambda+2\mu}{2}~\delta^{aA}.
\end{equation}
Similarly, for an initially stress-free material manifold, we
obtain
\begin{equation}
    \mathbbm{A}^{aA}{}_b{}^B=F^a{}_MF^c{}_Ng_{bc}\left[\lambda G^{AM}G^{BN}+\mu(G^{AB}G^{MN}+G^{AN}G^{BM})\right].
\end{equation}
For the case of an initially Euclidean material manifold with
Cartesian coordinates we have
\begin{equation}
    \mathbbm{A}^{aA}{}_b{}^B \frac{\partial^2 U^b}{\partial X^A\partial X^B}
    =(\lambda+\mu)U_{b,ab}+\mu U_{a,bb}.
\end{equation}
Therefore, Eq. (\ref{linearized-momentum}) reads
\begin{equation}
    (\lambda+\mu)U_{b,ab}+\mu U_{a,bb}=\frac{3\lambda+2\mu}{2}\frac{\partial \beta}{\partial x_a},
\end{equation}
where we have identified the indices $a$ and $A$. In analogy with thermal stresses, $\beta \delta_{ab}$ can be
thought of as an eigenstrain. See \citet{Goriely2008} for a review
of the existing linearized growth models.

\paragraph{Stress-free growth distributions in the linearized theory.} In this paragraph
we show that in dimension three if $\beta$ is linear in $\{X^A\}$,
i.e. if $\beta=\mathbf{a}\cdot\mathbf{X}$ for some constant vector
$\mathbf{a}$, then a stress-free body remains stress free after
growth. This is very similar to what is already known in classical linear
thermoelasticity: temperature distributions linear in Cartesian
coordinates leave a stress-free body stress free
\citep{BoleyWeiner1997,OzakinYavari2009}.

Let us consider a one-parameter family of material metrics
$\mathbf{G}_{\epsilon}$ and assume that the initial material
metric is Euclidean, i.e.
$\mathbf{G}_{\epsilon=0}=\boldsymbol{\delta}$. The corresponding
curvature tensor is $\boldsymbol{\mathcal{R}}_{\epsilon}$. We need to calculate the linearized curvature, i.e.
\begin{equation}
    \delta\boldsymbol{\mathcal{R}}=\frac{d}{d\epsilon}\Big|_{\epsilon=0}\boldsymbol{\mathcal{R}}_{\epsilon}.
\end{equation}
This will give the solution to stress-free growth distributions.
Note that $\delta\mathbf{G}=\frac{d}{d\epsilon}\big|_{\epsilon=0}\mathbf{G}_{\epsilon}$ corresponds to a linearized growth
and is stress-free if and only if $\delta\boldsymbol{\mathcal{R}}$
vanishes. To calculate the curvature variation, we follow
\citet{Hamilton1982} and denote derivative with respect to
$\epsilon$ by $'$. Following the definition of curvature tensor we
can write
\begin{equation}
    \mathcal{R}'_{ABCD}=-\frac{1}{2}\left(\frac{\partial^2 G'_{BD}}{\partial X^A \partial X^C}-\frac{\partial^2 G'_{BC}}{\partial X^A \partial X^D}
    -\frac{\partial^2 G'_{AD}}{\partial X^B \partial X^C}+\frac{\partial^2 G'_{AC}}{\partial X^B \partial X^D}\right)
    +\frac{1}{2}G^{PQ}\left(\mathcal{R}_{ABCP}G'_{QD}+\mathcal{R}_{ABPD}G'_{QC}\right).
\end{equation}
In the case of Ricci curvature
\begin{equation}
    R'_{AB}=G^{CD}\mathcal{R}'_{ACBD}+\left(G^{CD}\right)'\mathcal{R}_{ACBD}=G^{CD}\mathcal{R}'_{ACBD}-G^{CP}G^{DQ}G'_{PQ}\mathcal{R}_{ACBD}.
\end{equation}
Similarly, for scalar curvature we have
\begin{equation}
    \textsf{R}'=g^{AB}R'_{AB}+\left(G^{AB}\right)'R_{AB}=G^{AB}R'_{AB}-G^{AP}G^{BQ}G'_{PQ}R_{AB}.
\end{equation}

If the initial material manifold is Euclidean, i.e. if
$\mathcal{R}_{ACBD}=0$ and $R_{AB}=0$, we have
\begin{eqnarray}
  \delta\mathcal{R}_{ABCD} &=& -\frac{1}{2}\left(\frac{\partial^2 \delta G_{BD}}{\partial X^A \partial X^C}-\frac{\partial^2 \delta G_{BC}}{\partial X^A \partial X^D}
    -\frac{\partial^2 \delta G_{AD}}{\partial X^B \partial X^C}+\frac{\partial^2 \delta G_{AC}}{\partial X^B \partial X^D}\right), \\
  \delta R_{AB} &=& -\frac{1}{2}\left(\frac{\partial^2 \delta G_{CD}}{\partial X^A \partial X^B}-\frac{\partial^2 \delta G_{BC}}{\partial X^A \partial X^D}
  -\frac{\partial^2 \delta G_{AD}}{\partial X^B \partial X^C}+\frac{\partial^2 \delta G_{AB}}{\partial X^C \partial X^D}\right)\delta^{CD} , \\
  \delta \textsf{R} &=& \frac{\partial^2 \delta G_{BC}}{\partial X^A \partial X^D}\delta^{AB}\delta^{CD}-\frac{\partial^2 \delta G_{AB}}{\partial X^C \partial X^D}\delta^{AB}\delta^{CD}.
\end{eqnarray}

In the case of isotropic growth we have $\delta G_{AB}=\beta
\delta_{AB}$. In dimension three, vanishing of Ricci curvature is
equivalent to vanishing the curvature tensor. Thus, $\delta
R_{AB}=0$ reduces to
\begin{equation}
    \frac{\partial^2\beta}{\partial X^A\partial X^B}+\frac{\partial^2\beta}{\partial X^C\partial X^D}\delta^{CD}\delta_{AB}=0.
\end{equation}
This is equivalent to
\begin{eqnarray}
  && \beta_{,12}=\beta_{,13}=\beta_{,23}=0, \label{first-row}\\
  && 2\beta_{,11}+\beta_{,22}+\beta_{,33}=0, \label{row2} \\
  && \beta_{,11}+2\beta_{,22}+\beta_{,33}=0, \label{row3} \\
  && \beta_{,11}+\beta_{,22}+2\beta_{,33}=0 \label{row4}.
\end{eqnarray}
The three relations (\ref{first-row}) imply that
$\beta=f(X^1)+g(X^2)+h(X^3)$ for arbitrary functions $f,g$, and
$h$. The next three relations (\ref{row2})-(\ref{row4}) imply that
$\beta_{,11}=\beta_{,22}=\beta_{,33}=0$ and therefore
$f''(X^1)=g''(X^2)=h''(X^3)=0$, and hence $\beta$ is linear in
Cartesian coordinates of the initial material manifold.

In dimension two, $\delta\textsf{R}=0$ reduces to
\begin{equation}
    \frac{\partial^2\beta}{\partial X^A\partial X^B}\delta_{AB}=0.
\end{equation}
This means that $\beta$ has to be a harmonic function to represent
a stress-free growth distribution. Again, this is very similar to
what we know from classical linear thermoelasticity
\citep{BoleyWeiner1997,OzakinYavari2009}.

\section{Concluding Remarks}

In this paper, we presented a geometric theory of elastic solids with
bulk growth. We assumed that material points are preserved but
density and ``shape" are time dependent. We modeled a body with
bulk growth by a Riemannian material manifold with an evolving
metric tensor. The time dependency of material metric is such that
the growing body is always stress free in the material manifold.
We showed that energy balance needs to be modified when material
metric is time dependent. Covariance of energy
balance then gives all the balance laws. We also showed that entropy production
inequality has a non-standard form when material manifold has an
evolving metric. We showed that a more general notion of covariance of energy balance that includes temperature rescalings, in addition to giving all the balance laws, gives the constitutive restrictions imposed by the Clausius-Duhem inequality. We then showed how the principle of maximum entropy production can be used to obtain thermodynamically-consistent evolution equations for the material metric.

We showed how analytical solutions for the residual stress field
can be obtained in three examples of growing bodies with radial
symmetries. We showed that even if mass is conserved, i.e. when
growth results in only shape changes, still one may see residual
stresses. In the case of isotropic growth, we studied stress-free
growth distributions using the material curvature tensor in both
two and three dimensions.

A concrete connection was made between our geometric theory and the
conventional decomposition of deformation gradient into elastic
and growth parts. We showed that in a special coordinate basis
$\mathbf{F}_e$ is our $\mathbf{F}$. The present geometric theory
is more natural and does not introduce a mysterious intermediate
configuration. We linearized the nonlinear theory about a
reference motion. Assuming that both the ambient space and the
initial material manifold are Euclidean, we showed that growth
results in eigen strains very similar to those of classical linear
thermoelasticity. We found those growth distributions that are
stress free in the linearized framework in both dimensions two and three.

\section*{Acknowledgments}

The author benefited from discussions with A. Ozakin. He is also grateful to an anonymous reviewer for his/her good comments that led to the improvement of the original manuscript.


\end{document}